\documentclass[11pt]{article}
\usepackage[a4paper,left=3cm,top=3.5cm]{geometry}

\usepackage{amsmath,amscd,amssymb,amsfonts,mathdots,ntheorem}
\usepackage[colorlinks,linkcolor=blue,citecolor=red]{hyperref}
\usepackage{longtable}
\usepackage{cases}
\usepackage{tikz-cd}

\theoremstyle{plain}

\newtheorem{thm}{Theorem}[section] 

\newtheorem{defi}[thm]{Definition} 
\newtheorem{lemma}[thm]{Lemma} 

\newtheorem{cor}[thm]{Corollary}

\newtheorem{prop}[thm]{Proposition}
\newtheorem{rmk}[thm]{Remark}
\newtheorem*{ack}{Acknowledgement}
\newtheorem*{emp}{~~}
\newtheorem*{que}{Question}
\numberwithin{equation}{section}
\newtheorem{exa}[thm]{Example}
\def\comp{\ensuremath\mathop{\scalebox{.6}{$\circ$}}}
\def\mbc{\mathbf c}

\begin{document}
	


	
	\title{A Kodaira type conjecture on almost complex 4 manifolds}
	\date{}
	\author{Dexie Lin}
	\maketitle
	{College of Mathematics and Statistics, Chongqing University,
		Huxi Campus, Chongqing, 401331, P. R. China
		
		Chongqing Key Laboratory of Analytic Mathematics and Applications, Chongqing University, Huxi Campus, Chongqing, 401331, P. R.
		China}
	
 {E-mail: lindexie@126.com	}
	

	\begin{abstract}
		 Not long ago, Cirici and Wilson defined a Dolbeault cohomology on almost complex manifolds to answer Hirzebruch's problem.
		In this paper, we define  a refined Dolbeault cohomology on   almost complex manifolds. We show that the condition $\tilde h^{1,0}=\tilde h^{0,1}$ implies  a symplectic structure on a compact almost complex $4$ manifold, where $\tilde h^{1,0}$ and  $\tilde h^{0,1}$ are the dimensions of the refined Dolbeault cohomology groups with bi-degrees $(1,0)$ and $(0,1)$ respectively. Combining the partial answer to Donaldson's tameness conjecture, we offer a sufficient condition for a compact almost complex $4$ manifold to become  an  almost K\"ahler one. 
		 Moreover, we prove that the condition $\tilde{h}^{1,0}=\tilde h^{0,1}$ is equivalent to the generalized $\partial\bar\partial$-lemma. This can be regarded as an analogue of the Kodaira's conjecture on almost complex $4$ manifolds. As an application, we  show that the Kodaira-Thurston manifold satisfies  
		 the  $\partial\bar\partial$-lemma. Meanwhile, we show that   the Fr\"olicher-type equality does not hold on a  general almost complex $4$ manifold, which is  different to the case of  compact complex surfaces.  The main results of this paper can be regarded as the symplectic proof of the Kodaira's conjecture. 
	\end{abstract}
	\noindent
	{\bf Keywords}: almost complex manifold, Dolbeault cohomology, symplectic structure
	
	\noindent
	{\bf AMS classification}: 57R30, 53A45, 53C23, 53D35
	
	
	\section{Introduction}
	On   complex manifolds, the deRham differential has the decomposition	
	$d=\partial+\bar\partial,$
	where $\partial$ and $\bar\partial$ are first order differential operators with   bi-degrees $(1,0)$ and $(0,1)$ respectively, c.f., \cite[Definition 2.6.9 and Proposition 2.6.15]{Huy}. 
	The formula $\bar\partial^2=0$ allows a natural definition of the Dolbeault cohomology,   whose dimensions provide  complex invariants
	of  compact complex manifolds, c.f.,  \cite[Section 2 and 3-Chapter 0]{GH}. In this paper, a compact manifolds means a compact manifolds without boundary. 
	Similar to Riemannian manifolds,
	a  complex manifold $M$ equipped with a $J$-invariant metric is called an Hermitian manifold, where $J$ is the natural complex structure on $M$. On a compact Hermitian  manifold, the Hodge decomposition shows that $H^{*,*}_{Dol}$ are isomorphic to  ${\bar\partial}$-harmonic spaces with respect to the same bi-degrees, c.f., \cite[Section 6-Chapter 0]{GH}. 
	
	From the viewpoint of metrics,
	a K\"ahler manifold  is an Hermitian manifold of complex dimension $n$ such that for every point, there is a holomorphic coordinate chart around this point in which the metric agrees with the standard metric on $\mathbb C^n$ up to order $2$, c.f.,  \cite[Chapter 2 and Chapter 3]{Huy}. On   compact K\"ahler manifolds,   the K\"ahler identities(\cite[Proposition 3.1.12]{Huy}) and the Hodge decomposition yield
	\[ b_r =\sum_{j+k=r}h^{j,k}\mbox{ and }
	h^{j,k}=h^{k,j}.
	\]
	Here $h^{j,k}=\dim_{\mathbb C}H^{j,k}_{Dol}$ is the Hodge number of bi-degree $(j,k)$, c.f., \cite[Section 7-Chapter 0]{GH}. Consequently,  all
	Betti numbers with odd degrees   must be  even. 
	In \cite[Page 95]{KM},  Kodaira conjectured that: 
	\begin{emp}
		Every compact complex surface  admits a K\"ahler structure if and only if its first Betti number is even.
	\end{emp}
	By the Enrique-Kodaira's   classification \cite[Chapter VI]{BHPV} on compact complex surfaces,
	only the
	cases of elliptic surfaces and $K3$ surfaces remained open at that time.
	\begin{itemize}
		\item[(1)] Miyaoka \cite{Miy74}  proved that the
		conjecture holds for elliptic surfaces by using the current technique.   
		\item[(2)] Siu \cite{Siu}  overcame the difficulties in Todorov's proof \cite{Tod} and  proved 
		that every $K3$ surface is K\"ahler, hence completing the proof of Kodaira's	conjecture.
		\item[(3)] Buchdahl \cite{Buch} and Lamari \cite{Lamari} independently gave a unified proof of the Kodaira's conjecture.
	\end{itemize} 
	In this paper, we want to consider a variant of the  Kodaira's conjecture in the case of almost complex $4$ manifolds. Here, an almost complex manifold is  a manifold whose tangent bundle admits a complex structure.
	In a big picture, the Kodaira's  conjecture fills  in the condition for the upper-left arrow of the following diagram.
	
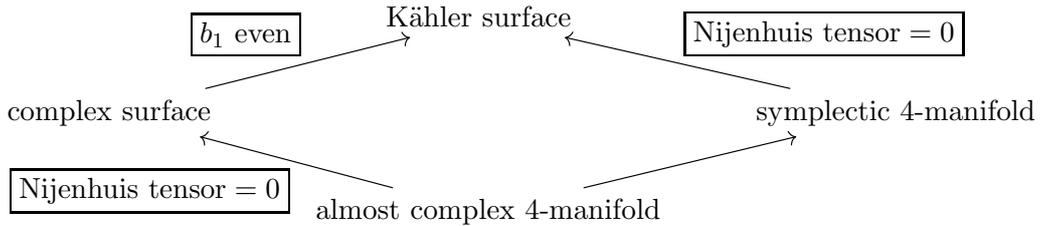
\[\begin{tikzcd}	    
		& \mbox{K\"ahler surface  }          &              \\
		\mbox{complex surface}\arrow[ru,"{\boxed{b_1\mbox{ even}} }"] &              & \mbox{symplectic 4-manifold} \arrow[ul,"\boxed{\mbox{Nijenhuis tensor}=0}"'] \\
		&\mbox{ almost complex 4-manifold}\arrow["\boxed{\mbox{Nijenhuis tensor}=0}",lu] \arrow[ru,""'] &  	 
	\end{tikzcd}\]	 
	
	From the above diagram,
	it is natural to ask: What is a proper condition for a compact almost complex $4$ manifold   to admit  a symplectic structure? Namely, what should be a proper generalization of the Kodaira's conjecture on compact almost complex $4$ manifolds? To start with,  the naive generalization is invalid. 
	In fact,   Gompf  \cite{Gom95} showed that: 
	\begin{emp}
		Every finitely presented group $G$ can be realized as  the fundamental group of some compact symplectic $4$ manifold.
	\end{emp} 
	By using Seiberg-Witten gauge theory, there are many  almost complex $4$ manifolds with  even first Betti numbers, but do not admit any symplectic structure, e.g., $(2k+1)\mathbb{C}P^2\# l\overline{\mathbb{C}P^2}$ with integers $k>0$ and $l\geq0$.
	Back to the case of  compact complex surfaces, the existence of a K\"ahler structure is equivalent to one of the following conditions,  c.f., \cite[Section 2 and 3-Chapter IV]{BHPV} and Gauduchon's work\cite{Gau76}:
	\begin{equation}
		(1)~b_1\mbox{ is even};~~~~(2)~\mbox{condition on }h^{1,0}\mbox{ and } h^{0,1};~~~~(3)~\partial\bar\partial-\mbox{lemma}.\label{eqn-kahler-condition}
	\end{equation}
	Here, the  $\mbox{condition on }h^{1,0}\mbox{ and } h^{0,1}$ means  one of the following equivalent equalities:
	\begin{equation}
		h^{1,0}=h^{0,1},~b_1=2h^{1,0}\mbox{ and }  b_1=2h^{0,1}.	\label{eqn-three-identities}
	\end{equation} 
Hence, it is natural to ask whether there is a suitable generalization of condition (2) or (3) in \eqref{eqn-kahler-condition} that serves as a sufficient condition for   the existence of a symplectic structure. 

To answer Hirzebruch's $20$-th. problem \cite{Hir54}, Cirici and Wilson \cite{CW21} defined a Dolbeault type cohomology for compact almost complex manifolds. However, Coelho,  Placini and Stelzig \cite{CPS} showed that the Hodge number $h^{0,1}$ is  finite if and only if such a compact almost complex $4$-manifold is a compact complex surface. In other words, the condition $h^{1,0}=h^{0,1}$ is equivalent to the K\"ahler condition on the compact almost complex $4$ manifolds.  Here $h^{1,0}$ and $h^{0,1}$ are the Hodge numbers of the  Dolbeault cohomology  defined by Cirici and Wilson.  
On the other hand, there are many compact symplectic $4$ manifolds without any K\"ahler structure  \cite{TO97}. Therefore, it raises a natural motivation  to expect another cohomology on compact almost complex manifold to study the symplectic structure. 
In this paper,  we define the  refined Dolbeault cohomology by  the cohomology of a complex $$\cdots\mathcal{A}^{*,*-1}_{Dol}\overset{\bar\partial}{\to}\mathcal{A}^{*,*}_{Dol}
\overset{\bar\partial}{\to}\mathcal{A}^{*,*+1}_{Dol}\cdots,$$
where $\mathcal{A}^{*,*}_{Dol}$ is a  certain subspace of the complex-valued forms, see Definition \ref{defi-tilde-Dolbeault}.  In particular, the identities $\partial^2|_{\mathcal{A}^{*,*}_{Dol}}=0=\bar\partial^2|_{\mathcal{A}^{*,*}_{Dol}}$ hold. This refined cohomology is denoted by $\tilde{H}^{*,*}_{Dol}$. 
 We give the refined Hodge numbers 
 on the Kodaira-Thurston manifold by the following diamond: 
\[\tilde{h}^{*,*}:=\dim_{\mathbb C}\Tilde{H}^{*,*}_{Dol}=\begin{array}{ccccc}
	& &1&& \\
	&\infty &&1&\\
	0&&\infty&&0\\
	&1&&1&\\
	&&1&&
\end{array}.\]
It is clear that this non-K\"ahlerian symplectic $4$ manifold
satisfies the identity $\Tilde{h}^{1,0}=\Tilde{h}^{0,1}$, see Example \ref{exa-kt} for more details. 
We start from the condition (2) in \eqref{eqn-kahler-condition} by  giving a sufficient condition for a compact almost complex $4$ manifold admitting a symplectic structure in terms of the refined Dolbeault cohomology. 
\begin{thm}\label{thm-main-1}
	Let $(X,J)$ be a compact almost complex $4$ manifold. Suppose that $\tilde h^{1,0}=\tilde h^{0,1}$. 
	Then, $X$ admits a $J$-taming  symplectic form. 
\end{thm}
Moreover, on compact almost complex $4$ manifolds, we show the following inequalities
\begin{equation}
	\tilde{h}^{1,0}=h^{1,0}\leq\tilde{h}^{0,1} \leq h^{0,1} ,	\label{eqn-tilde-Dolbeault-ineq}
\end{equation}
c.f.,    Lemma \ref{lemma-h1} and Corollary \ref{coro-betti-inequality}.

Next,   consider the  condition (3) in \eqref{eqn-kahler-condition}. 
Recall that on compact complex surfaces, Gauduchon \cite{Gau76} proved  that the  $\partial\bar\partial$-lemma is equivalent to any one of the three conditions in \eqref{eqn-three-identities}. Hence, the  $\partial\bar\partial$-lemma is equivalent to the existence of K\"ahler structure on compact complex surfaces. 
We  obtain an equivalent condition to the generalized $\partial\bar\partial$-lemma  on compact almost complex $4$ manifolds. 
\begin{thm}[generalized $\partial\bar\partial$-lemma]\label{thm-ddbar}
	Let $(X,J)$ be a compact almost complex $4$ manifold. Then,  any  $d$-exact $(1,1)$ form $\psi$ can be written as $\psi=\partial\bar\partial f$ for some    function $f$, if and only if the equality  $\tilde h^{1,0}=\tilde h^{0,1}$ holds.
\end{thm}
\noindent
Note that
	on compact K\"ahler manifolds, it is well-known that each $d$-exact form is $\partial\bar\partial$-exact, which is called the $\partial\bar\partial$-lemma.  Gauduchon's work \cite{Gau76} tells us that on the case of compact complex surfaces, the $\partial\bar\partial$-lemma is equivalent to the condition $h^{1,0}=h^{0,1}$. Hence, by the proof of Kodaira conjecture, the $\partial\bar\partial$-lemma is equivalent to the K\"ahler condition on  compact complex surfaces. 

The above results partially fill into the top arrow of the following  diagram.
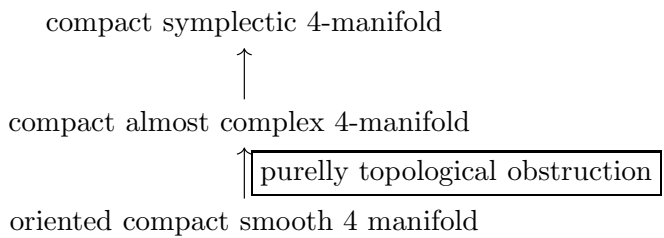
\begin{center}
	\begin{tikzcd}	    
		& \mbox{compact symplectic 4-manifold}        &              \\
		&    \mbox{compact almost complex 4-manifold  }   \arrow[u]    &  \\
		&\mbox{oriented compact smooth  4 manifold}\arrow[u,"\boxed{\mbox{purelly topological obstruction}}"',]&      
	\end{tikzcd}
\end{center}
Here, this topological obstruction is the existence of a certain cohomology class, c.f., {\cite[Proposition 9.3-Chapter:  IV]{BHPV}}.
For the  convenience of readers' understanding, we give the following diagram to illustrate the conditions for compact almost complex $4$ manifolds admitting the symplectic structure.

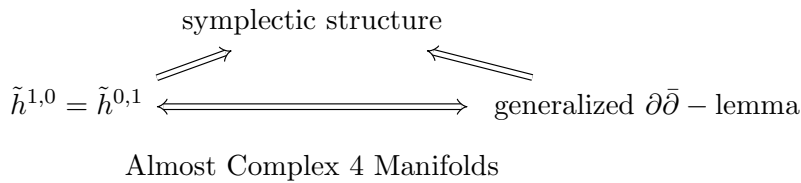
\begin{center}
	\begin{tikzcd}[column sep=0.02, row sep=small] 
		& {\begin{array}{l}
				\mbox{symplectic structure}
		\end{array}} &                                                   &                                                    \\
		{\tilde h^{1,0}=\tilde h^{0,1}} \arrow[rr, Leftrightarrow,"\begin{array}{c}
			~
			\\
			\mbox{Almost Complex 4 Manifolds}
		\end{array}"'] \arrow[ru, Rightarrow] &                           & \begin{array}{l}
		\mbox{generalized }	\partial\bar\partial- 
			\mbox{lemma}
		\end{array}  \arrow[lu, Rightarrow]                           
	\end{tikzcd}
\end{center} 
 Note that on compact complex surfaces, the three conditions are equivalent to each other, but for higher dimension there are non-K\"ahlerian complex manifolds satisfying the $\partial\bar\partial$-lemma, e.g., Moishezon manifolds.  Another motivation of this paper is to consider the question posted by Donaldson \cite{Don06}.  Hence, it is natural to consider the following question. 


\begin{que}
	Let $(X,\omega)$ be a compact symplectic $4$ manifold. Can we find a compactible almost complex structure $J$ on $X$,  such that  the equality
	\[
	\tilde h^{1,0}=\tilde{h}^{0,1},\]
	holds?
\end{que}


In the study   of $4$ manifolds, to determine whether a smooth manifold is symplectic or not is a popular topic. One can show that any non-compact almost complex $4$ manifold admits a symplectic structure by using Gromov's $h$-principle \cite[Page 84]{Gromov86}. Hence, the interesting topic is about the case of  compact $4$ manifolds.  There are topological constraints   for a certain kind of compact almost complex $4$ manifolds admitting a  symplectic structure by using the Seiberg-Witten gauge theory, see  Bauer's work \cite{Bau06} and Li's work \cite{Li06}. On the other hand, the exotic phenomenon  in smooth $4$ manifolds implies that there are also non-topological obstructions for the existence of a symplectic structure on a compact almost complex $4$ manifold, see   Fintushel and Stern's work \cite{FS98} on exotic $K3$
surfaces and Taubes' work \cite{Taub94}, \cite{Taub95}. 
A consequence of Theorem \ref{thm-main-1} and Theorem \ref{thm-ddbar}  is a sufficient condition for a compact  oriented smooth  $4$ manifold admitting a symplectic structure.

Similar to the case of compact complex manifolds, we give a Fr\"olicher type inequality on $b_1$, see Theorem \ref{thm-b1-inequality}. However, the equality can not hold in general, see the end of this paper.  
\noindent
In summary, we give a list of comparisons between the conditions for an almost complex $4$ manifold  admitting a  symplectic structure and the corresponding ones for a complex surface admitting a K\"ahler structure. 
\begin{center}
	\begin{tabular}{|c||c|}\hline  
		Compact    Complex  Surfaces& Compact Almost Complex $4$ Manifolds \\
		\hline
		\begin{tabular}{c}
			$b_1\equiv0\bmod2$  is equivalent  to  \\
			K\"ahler condition 
		\end{tabular}
		& \begin{tabular}{c}
			$b_1\equiv0\bmod2$  is unrelated to\\
			symplectic condition
		\end{tabular}   \\
		\hline
		\renewcommand{\arraystretch}{1.1}
		\begin{tabular}{c}
			$b_1=2h^{1,0}$   is equivalent  to  \\
			K\"ahler condition 
		\end{tabular}
		&\begin{tabular}{c}
			$b_1=2\tilde h^{1,0}=2h^{1,0}$  is unrelated to\\
			symplectic condition
		\end{tabular}\\\hline
		\renewcommand{\arraystretch}{1.1}
		\begin{tabular}{c}
			$h^{1,0}=h^{0,1}$  is equivalent  to  \\
			K\"ahler condition 
		\end{tabular}
		&\begin{tabular}{c}
			$\tilde h^{1,0}=\tilde h^{0,1}$  induces \\
			a symplectic structure
		\end{tabular} \\\hline
		\renewcommand{\arraystretch}{1.1}
		\begin{tabular}{c}
			$b_1=2h^{0,1}$  is equivalent  to  \\
			K\"ahler condition 
		\end{tabular} & \begin{tabular}{c}
			$b_1=\tilde h^{0,1}+\hat h^{0,1}$ induces\\
			a symplectic structure
		\end{tabular}
		\\\hline
		\renewcommand{\arraystretch}{1.1}
		\begin{tabular}{c}
			$\partial\bar\partial$-lemma  is equivalent  to  \\
			K\"ahler condition 
		\end{tabular} 
		&\begin{tabular}{c} generalized
			$\partial\bar\partial$-lemma induces  \\
			a symplectic structure
		\end{tabular}
		\\\hline
	\end{tabular}
\end{center}

\noindent
This paper is organized as follows: In Section 2, we review the related differential operators on almost complex manifolds, and the Dolbeault cohomology defined by Cirici and Wilson \cite{CW21}; In Section 3, we give the definition of  the refined Dolbeault cohomology and discuss some   relations between  the refined Dolbeault cohomology and the ones defined by Cirici and Wilson, we give an example to help readers' understanding; In Subsection 4.1, we show  some analogous  relations about the dimensions of the Dolbeault cohomology groups  on compact almost complex $4$ manifolds and prove the Fr\"olicher inequality; In Subsection 4.2, we give the proofs to Theorem \ref{thm-main-1} and Theorem \ref{thm-ddbar}.
\begin{ack}
	The author would like to express his gratitude to  Xiangdong Yang for introducing the Dolbeault cohomology  on almost complex manifold  and sharing   the work \cite{CW21} and \cite{ST23}, to Feng Wang for sharing \cite{Gom95} and Masaya Kwamura for explaining his work and Gauduchon's result \cite{Gau77}. The author also thanks Hongyu Wang for the helpful discussion and suggestion.  Part of the writing of this paper was completed at  Zhuhai Campus of Sun Yat-sen University, the author also thanks Yuan Wei for the hospitality. This paper is partially supported by NSFC No. 12301061. 
\end{ack}

\section{Operators on almost complex manifolds}

In this section, we give some necessary notions on almost complex  manifold, and introduce the related differential operators on  almost complex manifolds. We begin from the definition of almost complex manifold. 
\begin{defi}
	For a smooth manifold $M$ of real dimension $2n$. If there exists  a smooth section $J\in\Gamma(End(TM))$ on $TM$  satisfying $J^2=-1$, then we call
	$(M,J)$ an almost complex manifold with almost complex structure $J$.
\end{defi}
A celebrated theorem of  Newlander and Nirenburg  \cite{NN} states that an almost structure comes from a complex structure if and only if the associated Nijenhuis tensor 
\[N(v,w)=\frac{1}{4}\left([Jv,Jw]-[v,w]-J[v,Jw]-J[Jv,w]\right),\]
vanishes, where $v,w$ are vector fields on $M$.
Let  $T_{\mathbb C}M=TM\otimes_{\mathbb R}\mathbb C$ be the complexified tangent bundle. One has the decomposition with respect to $J$,
\[T_{\mathbb C}M=T^{1,0}M\oplus T^{0,1}M,\]
where $T^{1,0}M$ and $T^{0,1}M$ are the eigenspaces of $J$ corresponding to eigenvalues $i$ and $-i$ respectively.
The almost complex structure $J$ acts on the cotangent bundle $T^* M$ by $J  \alpha(v)=$ $-\alpha(J v)$, where $\alpha$ is a $1$-form and $v$ a vector field on $M$. This $J$  action  can be extended to any $p$-form $\psi$ by 
$$(J\psi)\left(v_1, \cdots, v_p\right)=(-1)^p \psi\left(J v_1, \cdots, J v_p\right).$$
Similarly, we  decompose $T^*_{\mathbb C}M:=T^{*}M\otimes_{\mathbb R}\mathbb C$ as
$T^*_{\mathbb C}M=T^{*,(1,0)}_{\mathbb C}M\oplus T^{*,(0,1)}_{\mathbb C}M$.
Denoting $ \mathcal A^*(M)$ by the space of sections of  $\bigwedge^*T^*_{\mathbb C}M$,   the decomposition holds:
\[ \mathcal A^r(M)=\bigoplus_{p+q=r} \mathcal A^{p,q}(M),\]
where $\mathcal{A}^{p,q}(M)$ denotes the space  of   forms of $\bigwedge^{p,q}T^*_{\mathbb C}M:=\bigwedge^pT^{*,(1,0)}_{\mathbb C}M\otimes\bigwedge^qT^{*,(0,1)}_{\mathbb C}M$.
By 
the formula \[d\psi\in   \mathcal{A}^{p+2,q-1}(M)+   \mathcal{A}^{p+1,q}(M)+
\mathcal{A}^{p,q+1}(M)+   \mathcal{A}^{p-1,q+2}(M),\]
for any $\psi\in   \mathcal{A}^{p,q}(M)$,
the deRham differential can be decomposed  as follows:
\[
d=\mu+\partial+\bar{\mu}+\bar{\partial} \mbox{ on } \mathcal{A}^{p,q},
\]
where   each  component is a derivation. The bi-degrees of the four components are  given by
\[|\mu|=(2,-1)
,~|\partial|=(1,0),~|\bar{\partial}|=(0,1),  ~|\bar{\mu}|=(-1,2).
\]
Observe that the operators $\partial,~\bar\partial$ are of the first order and the operators $\mu,~\bar\mu$   are of the zero order. 
Expanding  the relation $d^2=0$  gives the following: 
\begin{equation}
	\begin{aligned}
		\mu^{2} &=0, \\
		\mu \partial+\partial \mu &=0, \\
		\mu \bar{\partial}+\bar{\partial} \mu+\partial^{2} &=0, \\
		\mu \bar{\mu}+\partial \bar{\partial}+\bar{\partial} \partial+\bar{\mu} \mu &=0,\\
		\bar{\mu} \partial+\partial \bar{\mu}+\bar{\partial}^{2}&=0, \\
		\bar{\mu} \bar{\partial}+\bar{\partial} \bar{\mu}&=0, \\
		\bar{\mu}^{2}&=0.
	\end{aligned}\label{eqn-bar-operator}
\end{equation}
Similar to the case of Riemannian manifolds, we introduce the definitions of almost Hermitian manifold and almost K\"ahler manifold.
\begin{defi}
	A metric $g$   on an almost complex manifold $(M,J)$ is called    Hermitian, if $g$ is $J$-invariant, i.e., $g(J-,J-)=g(-,-)$.
	The imaginary component $\omega$ of $g$ defined by $\omega=g(J-,-)$, is a real non-degenerate $(1,1)$-form. 
	The quadruple   $(M,J,g,\omega)$ is called an \textit{almost Hermitian} manifold. 
	Moreover, when $\omega$ is $d$-closed, the quadruple   $(M,J,g,\omega)$ is called  an
	\textit{almost  K\"ahler} manifold. 
\end{defi}

We give some local calculation to help the understanding. 
Given an Hermitian metric $g$ on $(M,J)$, we set
$\{Z_r\}$ as a local $(1,0)$-frame with respect to  $g$ and  $\left\{\theta^r\right\}$ as the local associated coframe, i.e., $\theta^j(Z_k)=\delta^{j}_{k}$. Similar to the case of complex manifolds, we locally write $g_{k \bar j}=g\left(Z_k, \bar Z_j\right)$. The imaginary component $\omega$ of $g$   is  locally written as $\omega=i\sum_{k,j}g_{k \bar j}/2\theta^k \wedge \theta^{\bar{j}}$, c.f., \cite[Chapter 0-Section 6]{GH}. 
The coefficients of   Nijenhuis tensor $N$ is   given by
$$
N\left(Z_{\bar{j}}, Z_{\bar{k}}\right)=-\left[Z_{\bar{j}}, Z_k\right]^{(1,0)}= \sum_t N_{j k}^t Z_t,~ N\left(Z_j, Z_k\right)=-\left[Z_j, Z_k\right]^{(0,1)}=\sum_t\overline{N_{jk}^t} Z_{\bar{t}},
$$
and   the structure coefficients of the Lie bracket is expressed as $$[Z_j,\bar Z_{k}]=\sum_k C^r_{j\bar k}Z_r+
\sum_lC^{\bar l}_{i\bar j}Z_{\bar l},$$
c.f.,   \cite[2.8 and 2.9]{Kaw22}. Hence, we get the local expression of the Nijenhuis tensor
$$
N=\frac{1}{2}\sum_{j,k,l} \overline{N_{l j}^k} Z_{\bar{k}} \otimes\left(\theta^l \wedge \theta^j\right)+\frac{1}{2}\sum_{j,k,l} N_{l j}^k Z_k \otimes\left(\theta^{\bar{l}} \wedge \theta^{\bar{j}}\right) .
$$
The actions of the operators $\bar\mu,~\bar\partial$ and $\partial$ on $\theta^s$ take the following form:
$$\bar\mu \theta^s=\frac12\sum_{k,j}N^s_{\bar k\bar j}\bar\theta^k\wedge\bar\theta^j,~\bar\partial \theta^s=-\sum_{k,j}C^{s}_{k\bar j}\theta^k\wedge \bar\theta^j,~
\partial \theta^s=-\frac{1}{2}\sum_{k,j}C^s_{kj}\theta^k\wedge\theta^j.$$
By  conjugation, one has the similar expressions for $\mu\bar\theta^s,~\partial\theta^s$ and $\bar\partial\bar\theta^s$. By \cite[Lemma 2.1]{CW21}, it holds that 
\[
\mu+\bar\mu=N\otimes(Id_{\mathbb{C}}).\]	

The   Hermitian metric $g$ induces  a unique $\mathbb{C}$-linear Hodge-$*$ operator 
$$
*: \Lambda^{p, q} \rightarrow \Lambda^{n-q, n-p},
$$
$\mbox{defined by }
g\left(\phi_1, \phi_2\right) d V=\phi_1 \wedge *\overline{ \phi_2}$, 
where $dV$ is the volume form $\frac{\omega^n}{n!}$ and $\phi_1, \phi_2 \in \Lambda^{p, q}$. For any two forms $\varphi_1\mbox{ and }\varphi_2$,  we define their  product  by, 
$$(\varphi_1,\varphi_2)=\int_M\varphi_1\wedge *\overline{\varphi_2}.$$  
Let  $\delta^*$ denote the formal-$L^2$-adjoint of $\delta$ for $\delta\in\{\mu,\partial,\bar\partial,\bar\mu,d\}$ with respect to the metric $g$. 
By the integration by part, one has the identities 
\[\delta^*=-*\bar\delta*  \mbox{ for }\delta\in\{\mu,\partial,\bar\partial,\bar\mu,d\}.\] 
We also define the  Lefschetz operator
$$
L: \mathcal{A}^{p, q} \longrightarrow \mathcal{A}^{p+1, q+1},$$ 
by $ L(\alpha):=\omega \wedge \alpha$ for any form $\alpha\in\mathcal{A}^{p,q}$.
Its adjoint operator is $\Lambda=L^*=\star^{-1} L \star$. It is well known that the triple $\{L, \Lambda,  [L, \Lambda]\}$ defines a representation of $\mathfrak{s l}(2, \mathbb{C})$ and there is a Lefschetz decomposition on complex $k$-forms
$$
\mathcal{A}^k=
\bigoplus_{j \geq 0} 
L^jP^{k-2j},
$$
where $P^j=\ker(\Lambda) \cap \mathcal{A}^j$, c.f., \cite[Section 6, Chapter 0]{GH}. 

On almost  K\"ahler manifolds, Cirici and Wilson proved the  Hard Lefschetz Duality \cite[Theorem 5.1]{CW20} and the following two results. 

\begin{prop}[Cirici and Wilson {\cite[Proposition 3.1]{CW20}}]\label{prop-kahler-identities}
	On an   almost K\"ahler manifold $(M,J,g,\omega)$, i.e. $d\omega=0$, the following identities hold:
	\begin{itemize}
		\item[$(1)$]$[L, \bar{\mu}]=[L, \mu]=0$ and $\left[\Lambda, \bar{\mu}^*\right]=\left[\Lambda, \mu^*\right]=0$.
		\item[$(2)$] $[L, \bar{\partial}]=[L, \partial]=0$ and $\left[\Lambda, \bar{\partial}^*\right]=\left[\Lambda, \partial^*\right]=0$.
		\item[$(3)$]  $\left[L, \bar{\mu}^*\right]=i \mu,\left[L, \mu^*\right]=-i \bar{\mu}$ and $[\Lambda, \bar{\mu}]=i \mu^*,[\Lambda, \mu]=-i \bar{\mu}^*$.
		\item[$(4)$]  $\left[L, \bar{\partial}^*\right]=-i \partial,\left[L, \partial^*\right]=i \bar{\partial}$ and $[\Lambda, \bar{\partial}]=-i \partial^*,[\Lambda, \partial]=i \bar{\partial}^*$.
	\end{itemize}
\end{prop}


\begin{thm}[Cirici and Wilson {\cite[Theorem 4.1]{CW20}}]\label{thm-CW-02}
	On any compact almost K\"ahler manifold with dimension $2 n$, it holds
	$$
	\mathcal{H}_d^{p, q}=\mathcal{H}_{\bar\partial}^{p, q} \cap \mathcal{H}_\mu^{p, q}=\mathcal{H}_{\partial}^{p, q} \cap \mathcal{H}_{\bar{\mu}}^{p, q},
	$$
	for all $(p, q)$. Moreover,
	we also have the following identities:
	\begin{itemize}
		\item[$(1)$](Complex conjugation)
		$$
		\mathcal{H}_{\bar{\partial}}^{p, q} \cap \mathcal{H}_\mu^{p, q}=\mathcal{H}_{\bar{\partial}}^{q, p} \cap \mathcal{H}_{{\mu}}^{q, p} .
		$$
		\item[$(2)$]  (Hodge duality) 
		$$
		*: \mathcal{H}_{\bar{\partial}}^{p, q} \cap \mathcal{H}_\mu^{p, q} \rightarrow \mathcal{H}_{\bar{\partial}}^{n-q, n-p} \cap \mathcal{H}_\mu^{n-q, n-p} .
		$$
		\item[$(3)$] (Serre duality)
		$$
		\mathcal{H}_{\bar\partial}^{p, q} \cap \mathcal{H}_\mu^{p, q} \cong \mathcal{H}_{\bar\partial}^{n-p, n-q} \cap \mathcal{H}_\mu^{n-p, n-q} .
		$$
	\end{itemize}
\end{thm}

For the convenience of  later arguments in the next section, we set $$\mathcal{H}^{p,q}_{\delta,\delta'}:=\mathcal{H}^{p,q}_\delta\cap\mathcal{H}^{p,q}_{\delta'},$$
for $\delta,\delta'\in\{\partial,\bar\partial,\mu,\bar\mu\}$ and $\ell^{p,q}=\dim_{\mathbb C} \mathcal{H}^{p,q}_{\bar\partial,\mu}$. 

\section{Cohomology groups on almost complex manifolds}
First, we review the Dolbeault cohomology $H^{*,*}_{Dol}$ defined by Cirici and Wilson on   almost complex manifolds. 

Recall that 
the  last three identities in \eqref{eqn-bar-operator} imply  the following identities:
\[
\bar{\mu} \bar{\partial}+\bar{\partial} \bar{\mu}=0, ~
\bar{\mu}^{2}=0, \mbox{ and }\bar \partial^2|_{\ker(\bar\mu)}\equiv0(\bmod ~ im(\bar\mu)). \]
These identities induce a  Fr\"olicher-type spectral sequence on $\mathcal{A}^{*,*}$. 
Hence, one can define the Dolbeault cohomology for an  almost  complex manifold
$(M,J)$.   
 If there is no ambiguity, we omit the underline manifold and almost complex structure and  simply write  the cohomology as $H^{*,*}_{Dol}$
\begin{defi}[Cirici and Wilson {\cite[Definition 3.1]{CW21}}]
	The Dolbeault cohomology of an almost  complex $2n$-dimensional manifold
	$(M,J)$ is given  by
	\[
	H^{p,q}_{Dol}=H^q(H^{p,*}_{\bar\mu},\bar\partial)=
	\frac{\ker(\bar\partial:H^{p,q}_{\bar\mu}\to H^{p,q+1}_{\bar\mu})}
	{im(\bar\partial:H^{p,q-1}_{\bar\mu}\to H^{p,q}_{\bar\mu})},
	\]
	where $H^{p,q}_{\bar\mu}=\frac{\ker(\bar\mu:\mathcal{A}^{p,q}\to \mathcal{A}^{p-1,q+2})}{im(\bar\mu:\mathcal{A}^{p+1,q-2}\to \mathcal{A}^{p,q})}$.
\end{defi}
Considering the presence of the operator $\bar\mu$,
Cirici and Wilson gave the definition of the Hodge-type filtration.

\begin{defi}[Cirici and Wilson {\cite[Definition 3.2]{CW21}}]
	The Hodge filtration of $\mathcal{A}^{*}$ of an almost complex manifold $(M,J)$ is given by the following decreasing filtration 
	$$
	F^p \mathcal{A}^n:=\operatorname{Ker}(\bar{\mu}) \cap \mathcal{A}^{p, n-p} \oplus \bigoplus_{i>p} \mathcal{A}^{i, n-i}.
	$$
\end{defi}
Cirici and Wilson  showed that the Dolbeault cohomology groups are isomorphic to the first page of the above Hodge filtration, and the spaces of  $(\Delta_{\bar\partial}+\Delta_{\bar\mu})$-harmonic forms  can be embedded into  the Dolbeault cohomology groups with respect to the same bi-degree. 

\begin{thm}[Cirici and Wilson  {\cite[Theorem 3.8 and Proposition 4.10]{CW21}}]\label{thm-CW-21}
	Let $(M,J)$ be  a $2n$-dimensional compact almost complex  manifold. One has the following:
	\begin{itemize}
		\item[$(1)$] \[
		H^{p,q}_{Dol}\cong E^{p,q}_1:=
		\frac{\{\omega\in\ker(\bar\mu)\cap\mathcal{A}^{p,q}|~\bar\partial\omega\in im(\bar\mu)\}}
		{\{im({\bar\partial})  \cap\ker(\bar\mu)+im(\bar\mu)\}\cap\mathcal{A}^{p,q}}\Rightarrow H^{p+q}_{dR}, ~0\leq p,q\leq n,
		\]
		where $E^{p,q}_1$ is the first stage of  the spectral sequence of the above Hodge filtration;
		\item[$(2)$] \[
		\ker(\Delta_{\bar\partial}+\Delta_{\bar\mu})\cap\mathcal{A}^{p,q}\subseteq H^{p,q}_{Dol} \mbox{ for any }p,~q,
		\]
		and the equality holds for $p\in\{0,...,n\}$ and $q\in\{0,n\}$.
	\end{itemize}
\end{thm}

Inspired by the work of Sillari and  Tomassini \cite{ST23}, we define the certain subspaces of $\mathcal{A}^{*,*}$ and the associated Dolbeault type cohomology. Here, we also omit the underline manifold and the almost complex, just write $\tilde{H}^{*,*}_{Dol}$. 
\begin{defi}\label{defi-tilde-Dolbeault}
	On an almost complex manifold $(M,J)$,	we define the subspace $\mathcal{A}^{*,*}_{Dol}$ of $\mathcal{A}^{*,*}$ by
	\[
	\mathcal{A}^{*,*}_{Dol}=\ker(\mu)\cap \ker(\bar\mu)
	\cap \ker((\bar\partial+\mu)\bar\partial)\cap
	\mathcal{A}^{*,*}.
	\]
	We define the refined Dolbeault cohomology  $\tilde{H}^{*,*}_{Dol}$ by
	\[\tilde{H}^{*,*}_{Dol}=\frac{\ker(\bar\partial)\cap\mathcal{A}^{*,*}_{Dol}}{\bar\partial(\mathcal{A}^{*,*-1}_{Dol})}.\]
\end{defi}
One can check that the identities $\partial^2|_{\mathcal{A}^{*,*}_{Dol}}=0\mbox{ and }
\bar\partial^2|_{\mathcal{A}^{*,*}_{Dol}}=0$ hold, and $\tilde{H}^{*,*}_{Dol}$ is the cohomology of the complex sequence 
\begin{equation}
	\cdots\to \mathcal{A}^{p,q-1}_{Dol}\overset{\bar\partial}{\to}
	\mathcal{A}^{p,q}_{Dol}\overset{\bar\partial}{\to}\mathcal{A}^{p,q+1}_{Dol}\to\cdots. \label{eqn-tilde-complex}
\end{equation}

\noindent
By observation, we have the following relations on compact almost Hermitian manifolds. 
\begin{lemma}   On a  compact almost Hermitian manifold $(M,J,g,\omega)$ with dimension $2n$, we have the following relations:
	\begin{itemize}
		\item[$(1)$] $\tilde H^{p,0}_{Dol}\cong H^{p,0}_{Dol}\cong \mathcal{H}^{p,0}_{\bar\partial,\bar\mu}$, for $0\leq p\leq n$.
		\item[$(2)$] $\mathcal{H}^{0,1}_{\bar\partial,\mu}\subseteq \tilde H^{0,1}_{Dol} \subseteq H^{0,1}_{Dol}$; 
	\end{itemize}
\end{lemma}
\begin{pf}  
	By the identification $$H^{p,q}_{Dol}\cong     \frac{\{\sigma\in\ker(\bar\mu)\cap\mathcal{A}^{p,q}|~\bar\partial\sigma\in im(\bar\mu)\}}    {\{im({\bar\partial})  \cap\ker(\bar\mu)+im(\bar\mu)\}\cap\mathcal{A}^{p,q}},$$
	the linear map $\alpha\in \ker(\mu)\cap\ker(\bar\mu)\cap \ker((\bar\partial+\mu)\bar\partial)\cap\ker(\bar\partial)\cap\mathcal{A}^{p,q} \mapsto [\alpha] \in H^{p,q}_{Dol}$   gives a well-defined  linear map $\tilde H^{p,q}_{Dol}\to H^{p,q}_{Dol}$ for each  $(p,q)$. 
	\begin{itemize}
		\item[$(1) $] The isomorphism $\tilde H^{p,0}_{Dol}\cong H^{p,0}_{Dol}$ follows from 
		the facts that the operators $\mu$ and $\bar\mu^*$ act  trivially on $\mathcal{A}^{p,0}$ and the intersection $im(\bar\partial)\cap\mathcal{A}^{p,0}$ is trivial.
		\item[$(2)$] When a class $[\alpha]\in\tilde H^{0,1}_{Dol}$ is trivial in $H^{0,1}_{Dol}$,  we write $\alpha=\bar\partial f$ for some function $f$. By the definition of $\alpha$,  one gets that $f\in\mathcal{A}^0_{Dol}$, i.e., $[\alpha]$ is trivial in $\tilde H^{0,1}_{Dol}$. Similarly, if a form $\alpha\in\mathcal{H}^{0,1}_{\bar\partial,\mu}$ is trivial in $\tilde H^{0,1}_{Dol}$, then $\alpha=\bar\partial f$.
		Applying the maximum principle to the equation $\bar\partial^*\bar\partial f=0$,   $f$ is a constant function and $\alpha=0$.
	\end{itemize}	
\end{pf}

\noindent
To compare the dimensions of  Dolbeault cohomology groups with the first Betti number,   we introduce the following  two terms.
\begin{defi}\label{defi-hat-Dolbeault}
	We define two terms $\hat{H}^1$ and $\hat{H}^{0,1}$ by
	\[\hat{H}^1=\frac{\{u_1+u_2\in\mathcal{A}^{1,0}\oplus\mathcal{A}^{0,1}\mid~
		\bar\mu(u_1)+\bar\partial(u_2) =0=\mu (u_2)+\partial (u_1)\}}
	{\{u'_1+u'_2\in( \partial(\mathcal{A}^0)\oplus \bar\partial(\mathcal{A}^0))\mid
		\bar\mu(u'_1)+\bar\partial(u'_2)=\mu (u'_2)+\partial (u'_1) =0	\}},\]
	and  \[\hat H^{0,1}=\frac{\{u''\in\mathcal{A}^{0,1}|~d^{2,0}(u'+u'')=0=d^{0,2}(u'+u'')\mbox{ for some }u'\in\mathcal{A}^{1,0}\}}{im(\bar\partial)\cap\mathcal{A}^{0,1}}.\]
\end{defi} 
By  definition, it is clear that $\hat H^{0,1}\subseteq H^{0,1}_{Dol}$. For two functions $f_1$ and $f_2$, the condition $\bar\mu(\partial f_1)+\bar\partial(\partial f_2) =\mu (\partial f_2)+\partial (\partial f_1)=0$ is equivalent to the condition $\partial^2(f_1-f_2)=\bar\partial^2(f_1-f_2)=0$, which  trivially holds on complex manifolds. 
We denote their dimensions by 
\[\hat h^{0,1}=\dim_{\mathbb C} \hat H^{0,1},~\hat h^1=\dim_{\mathbb C}\hat H^1, \tilde h^{p,q}=\dim_{\mathbb C}\tilde H^{p,q}_{Dol} \mbox{ and }
h^{p,q}=\dim_{\mathbb C}  H^{p,q}_{Dol}.
\] 	

Now, we give a relation between $H^1_{dR}$ and $\hat{H}^1$ on compact almost complex manifolds. 
\begin{lemma}
	On a compact almost complex manifold, there is an inclusion \[
	H^{1}_{dR}\hookrightarrow \hat{H}^1.\]
\end{lemma}
\begin{pf}
	Notice that the composition of the  inclusion and the projection $$\ker(d)\cap(\mathcal{A}^{1,0}\oplus\mathcal{A}^{0,1})\to\ker(d^{2,0})\cap \ker(d^{0,2})\cap(\mathcal{A}^{1,0}\oplus\mathcal{A}^{0,1})\to\mathcal{A}^{0,1},$$ 
	gives a canonical map $H^{1}_{dR}\to \hat{H}^1$, where $d^{2,0}$ and $d^{0,2}$ are  $(2,0)$ and $(0,2)$-components of $d$ respectively. 
	
	It suffices to show that this map is an injection. 
	If  a class $[\alpha]\in H^{1}_{dR}$ maps to a trivial class in $\hat{H}^1$. Then, $\alpha=\partial f_1+\bar\partial f_2$ for   two functions $f_1$ and $f_2$, where $\alpha$ is a representative of this class. Since $d\alpha=0$, one gets that 
	\[\bar\partial\partial f_1+\partial\bar\partial f_2=0, \]
	hence $\partial\bar\partial(f_1-f_2)=0$ by the formula $\partial\bar\partial+\bar\partial\partial=0$ on functions. Applying the maximum principle to the equation $i\partial\bar\partial(f_1-f_2)=0$,  one obtains that $f_1-f_2$ is a constant. Therefore, we write $\alpha=\partial f_1+\bar\partial f_2=\partial f_1+\bar\partial f_1=df_1$, which is trivial in $H^{1}_{dR}$.
\end{pf}

For the convenience of the later application in the next section, we end this section by the following lemma.

\begin{lemma}\label{lemma-b1-h1}
	The  equality
	\begin{equation} \hat h^{1}= \hat h^{0,1}+\tilde h^{0,1},\label{eqn-b1-h1}
	\end{equation}	holds
	on any compact almost complex manifold.
\end{lemma}
\begin{pf}
	First, we define a linear map 
	\[
	\ker(d^{0,2})\cap\ker(d^{2,0})\cap\mathcal{A}^{1}\to 
	\mathcal{A}^{0,1},
	\]
	by $u'+u''\mapsto u''$,  where $u'\in\mathcal{A}^{1,0}$ and  $u''\in\mathcal{A}^{0,1}$. 
	Indeed, the above map induces a well-defined surjective linear map $$\pi:\hat H^1\to \hat H^{0,1},$$ 
	$\mbox{defined by }[u'+u'']\mapsto[u'']$. 
	
	If a class $[u'+u'']$ belongs to the kernel of $\pi$,  then
	$u''=\bar\partial f$ for some function $f$.  
	Since $[u'+\bar\partial f]\in\hat H^{1}$, one gets 
	\[ \bar\partial^2 f+\bar\mu u'=0\mbox{ and }\partial u'+\mu\bar\partial f=0.
	\]
	Substituting  the identities $\partial^2+\mu\bar\partial=0$ and $\bar\partial^2+\bar\mu\partial=0$ on functions into the above two equations gives that 
	\[u'-\partial f\in\ker(\bar\mu)\cap\ker(\partial),
	\]
	i.e., $[{\bar u'-\bar \partial \bar f}]\in \tilde H^{0,1}_{Dol}$.
	It is clear that $\overline{\tilde H^{0,1}_{Dol}}=\tilde{H}^{1,0}_{\partial}:=\frac{\{v\in\mathcal{A}^{1,0}\mid \partial(v)=0=\bar\mu(v)\}}{ \partial(\mathcal{A}^0_{D0l})\cap\mathcal{A}^{1,0}_{Dol}}$. 	On the other hand, we define a map
	\[
	\tilde{H}^{1,0}_{\partial}\to \hat H^1,\]
	$\mbox{by }
	[v]\mapsto[v+0]$. This map is an injection, because if a class  $[v]\in\tilde{H}^{1,0}_{\partial}$ maps to a trivial class in $\hat{H}^1$, then $v=\partial f$ with $\partial f\in\ker(\partial)\cap\ker(\bar\mu)$ by $\partial(v)=\bar\mu(v)=0$, which is trivial in $\tilde{H}^{1,0}_{\partial}$.
	Hence, for any class $[v]\in {\tilde H^{0,1}_{\partial}}$, the class $[ v+ u'+\bar\partial  f]$ belongs to $\ker(\pi)$. Therefore, we get the exact sequence
	\[0\to\overline{\tilde{H}^{0,1}_{Dol}}=\tilde{H}^{1,0}_{\partial}\to \hat H^1\to \hat H^{0,1}\to 0,
	\]
	i.e., $\hat h^{1}=\hat h^{0,1}+\tilde h^{0,1}$. 
\end{pf}



Recall that on a compact almost complex $4$ manifold with a non-trivial almost complex structure, $H^{2,0}_{Dol} = H^{0,2}_{Dol}=0$ by Cirici and Wilson \cite{CW20}. 
We show  $\tilde{H}^{1,0}_{Dol}$ and $\tilde{H}^{0,1}_{Dol}$ are  trivial on  some compact almost complex manifolds with $\dim_{\mathbb R}\geq6$. 

\begin{lemma}\label{lemma-trivial}
	Let $(M,J)$ be a   compact   almost complex  $2n(\geq6)$-manifold   such that the associated  Nijenhuis tensor is not vanishing and takes the maximal rank at some point. Then, $\tilde{H}^{1,0}_{Dol}=\{0\}$. Moreover, if the Nijenhuis tensor takes the maximal rank anywhere, then $\tilde{H}^{0,1}_{Dol}=\{0\}$.  
\end{lemma}
\begin{pf}
	An easy   calculus yields $$rank_{\mathbb C}\mathcal{A}^{1,0}=rank_{\mathbb C}\mathcal{A}^{0,1}=n,~rank_{\mathbb C}\mathcal{A}^{2,0}=rank_{\mathbb C}\mathcal{A}^{0,2}=\frac{n(n-1)}{2}.$$ 
	For an Hermitian metric $g$, it is clear that $\bar\partial^* u=0$ and $u\in\ker(\bar\partial)$ is equivalent to $u\in\ker(\Delta_{\bar{\partial}})$ for any $u\in \mathcal{A}^{1,0}$. 
	Hence, all elements in $\tilde{H}^{1,0}$ are $\Delta_{\bar{\partial}}$-harmonic. Since $N\not\equiv0$, there exists a small open subset $U$ such that the restriction $N|_U$ is nowhere vanishing. 
	When $n\geq3$, it is clear that $\bar\mu|_U:\mathcal{A}^{1,0}(U)\to\mathcal{A}^{0,2}(U)$ is either  isomorphic($n=3$) or injective($n>3$) by the assumption on the maximal rank. The equation $\bar\mu(u)=0$ gives that $u|_U=0$ for any $u\in\mathcal{A}^{1,0}$.
	Hence,    the unique-continuity-property for harmonic forms shows that $u=0$ for any $u\in\tilde{H}^{1,0}_{Dol}$. 
	
	Similarly, 
	when $N$ is nontrivial and takes the maximal rank everywhere,  $\mu:\mathcal{A}^{0,1}\to\mathcal{A}^{2,0}$ is either  isomorphic($n=3$) or injective($n>3$). Therefore, $\mu(v)=0$ implies $v=0$ for any $v\in\mathcal{A}^{0,1}$, which proves the last statement. 
\end{pf}
\def\bp{\bar\theta}

We end this section by giving an explicit calculation of a compact almost complex  manifold. 

\def\mbc{\mathbb{C}}
\begin{exa}\label{exa-kt}
	The Kodaira-Thurston manifold $\mathrm{KT}^4$ is defined to be the  product manifold $S^1 \times  H_3(\mathbb{Z}) \backslash H_3(\mathbb R) $, where $H_3(\mathbb{R})\subset G L(3, \mathbb{R})$ denotes the Heisenberg group 
	and $H_3(\mathbb{Z})=H_3(\mathbb{R}) \cap G L(3, \mathbb{Z})$ acts on $H_3(\mathbb{R})$ by left multiplication. 
	This manifold 
	can be  given by identifying points in $\mathbb R^4$ with the relation
	$$(t,x,y,z)\sim (t+t_0 ,
	x+x_0,
	y+y_0,
	z+z_0+x_0 y),
	$$
	for each element $(t_0, x_0, y_0, z_0) \in \mathbb{Z}^4$ and $(t,x,y,z)\in\mathbb R^4$.
	Choose an almost complex structure given by
	\[J\frac{\partial}{\partial t}=\frac{\partial}{\partial x},~J\frac{\partial}{\partial x}=-\frac{\partial}{\partial t},~J(\frac{\partial}{\partial y}+x\frac{\partial}{\partial z})=\frac{\partial}{\partial z},~J\frac{\partial}{\partial z}=-(\frac{\partial}{\partial y}+x\frac{\partial}{\partial z}).\]
	The vector fields$$v_1=\frac{1}{2}\left(\frac{\partial}{\partial t}-i \frac{\partial}{\partial x}\right)\mbox{ and } v_2=\frac{1}{2}\left(\left(\frac{\partial}{\partial y}+x \frac{\partial}{\partial z}\right)-i\frac{\partial}{\partial z}\right),$$ span  $T^{1,0} M$ at each point. 
	The  $1$-forms
	$$
	\theta_1=dt+i dx \mbox{ and }\theta_2=dy+i (dz-x dy) ,
	$$ are dual to $v_1$ and $v_2$ respectively. It is clear that 
	$\theta_1$ and $\theta_2$ satisfy the structure equations
	$$
	\begin{gathered}
		d \theta_1=0,
		\partial \theta_2=-\frac{1}{4} (\theta_1 \wedge \theta_2),~\bar\partial\theta_2=-\frac{1}{4}(\theta_1 \wedge \bar{\theta}_2+\theta_2 \wedge \bar{\theta}_1),~\bar\mu\theta_2=\frac{1}{4}(\bar{\theta}_1 \wedge \bar{\theta}_2 ) .
	\end{gathered}
	$$
	$J$ has a compactible symplectic form $\omega=i/2(\theta_1\wedge\bar\theta_1+\theta
	_2\wedge\bar \theta_2)$, i.e. $(KT^4, g_J,J,\omega)$ is an almost K\"ahler $4$ manifold, where $g_J=\omega(J,)$. It is known that this symplectic manifold has no K\"ahlerian structure by $b_1=3$. 
	
	We start from $\tilde{H}^{0,1}_{Dol}$ to compute all  refined Hodge numbers. 
	For any $(0,1)$-form $u=f_1\bp_1+f_2\bp_2$ with smooth functions $f_1$ and $f_2$, the requirement $\mu(u)=0=\bar\partial(u)$ implies that $$f_2=0,~\bar\partial f_1\wedge\bp_1=-\frac{1}{2}\left(\left(\partial_y+x  \partial_z\right)f_1+i  \partial_zf_1\right)\bp_1\wedge\bp_2=0,$$
	where $\partial_y=\frac{\partial}{\partial y}$ and $\partial_z=\frac{\partial}{\partial z}$.
	On the other hand, for any smooth function $f$,we write $\partial f=f_1\theta_1+f_2\theta_2$, where $f_i=v_i(f)$ for $i=1,2$. The requirement $\bar\partial^2 f=0$ derives
	$v_2(f)=0$, and the requirement $\partial^2 f=0$ derives $\bar v_2f=0$. 
	We regard $f$ as a smooth bounded function on $\mathbb R^4$, still denoted by $f$. Setting $$f'(t,x,y,z)=f(t,x,y,xy+z),$$ the equation $\bar v_2f=0$ derives 
	\[\partial_y f'+i\partial_zf'=0. \]
	Namely, the function
	$f'_{t,x}(y,z):=f'(t,x,y,z)$ is  holomorphic and bounded for any $(t,x)$. This implies that $f'$ is independent of the coordinates $\{y,z\}$ by the Liouville theorem. In particular, it holds that $\partial_yf'=\partial_zf'=0$, i.e. $$\partial_y f=\partial_zf=0.$$ 
	Thus, the equation $\bar v_2f=0$ implies that $f$  depends only on  $\{t,x\}$.
	In fact, one can use the Fourier expansion(c.f. \cite[Proposition 3.1 and 3.2]{HZ20}) on $KT^4$ to obtain the same result.  
	
	
	Therefore,
	we reduce the calculation of the cohomology to the following quotient on the torus $T^2$, 
	\[\Tilde{H}^{0,1}_{Dol}\cong \frac{\{C^\infty(T^2)\}}{\{im(v_1)\}}.\]
	The right-hand-side part of the above identity has the  inclusions 
	\[\mbc [\bar\theta_1]\subseteq\frac{\{C^\infty(T^2)\}}{\{im(v_1)\}}\subseteq \frac{\{C^\infty(T^2)\}}{im(v_1\bar v_1)}.\]
	Since the operator $v_1\bar v_1=\frac{\partial^2}{\partial t^2}+\frac{\partial^2}{\partial x^2}$ is a formal self-adjoint elliptic operator on the torus, its cokernel is of $\dim=1$ by its index and the maximum-principle. Hence, the relations $0<\Tilde{h}^{1,0}\leq b_1/2$ and $\dim_\mbc\Tilde{H}^{0,1}_{Dol}\geq \Tilde{h}^{1,0}$ imply that $ \Tilde{h}^{0,1}=\tilde{h}^{1,0}=1$.
	
	To calculate $\Tilde{H}^{1,1}_{Dol}$, we write 
	\[w=a_1\theta_1\wedge\bar\theta_1+a_2\theta_2\wedge\bar\theta_2+b_1\theta_1\wedge \bar\theta_2+b_2\theta_2\wedge\bar\theta_1,\]for any $(1,1)$-form $w$.
	The equation $\bar\partial w=0$ gives 
	\[\bar v_2a_1\bar\theta_2\wedge \theta_1\wedge\bar\theta_1+\bar v_1a_2\bar\theta_1\wedge\theta_2\wedge \bar\theta_2+
	\bar v_1b_1 \bar \theta_1\wedge \theta_1\wedge \bar\theta_2-\frac{b_1}{4}\theta_1\wedge\bar\theta_1\wedge\bar\theta_2+\bar v_2b_2\bar\theta_2\wedge\theta_2\wedge\bar\theta_1-\frac{b_2}{4}\theta_1\wedge\bar\theta_2\wedge\bar\theta_1=0,\]
	i.e. \[\begin{cases}
		\bar v_2a_1-\bar v_1b_1-\frac{b_1}{4}+\frac{b_2}{4}=0,\\
		-\bar v_1a_2+\bar v_2b_2=0.
	\end{cases}\]
	Consider   $u=f_1\theta_1+f_2\theta_2\in\mathcal{A}^{1,0}\cap\ker(\bar\mu)\cap\ker(\bar\partial^2)$ for some functions $f_1$ and $f_2$. By the similar arguments, one has $f_2=0$ and 
	$ v_2f_1=0$, i.e.,  $f_1$ is a function depending only on $\{t,x\}$. 
	This establishes the following inclusion
	\[
	\tilde{H}^{1,1}_{Dol}\supset
	\frac{\{a(t,x)\theta_1\wedge\bar\theta_1+b_1(t,x)\theta_1\wedge\bar\theta_2+b_2(t,x)\theta_2\wedge\bar\theta_1\mid ~b_2=b_1 +4\bar v_1b_1 \}}{\{im(\bar  v_1(C^\infty(T^2)))\theta_1\wedge\bar\theta_1\}},
	\]
	where the functions $a(t,x)$,  $b_1(t,x)$ and $b_2(t,x)$ depend only on $\{t,x\}$.
	It is clear that \[
	\dim_\mbc\frac{\{a(t,x)\theta_1\wedge\bar\theta_1+b_1(t,x)\theta_1\wedge\bar\theta_2+b_2(t,x)\theta_2\wedge\bar\theta_1\mid ~b_2=b_1+4\bar v_1b_1 \}}{\{im(\bar  v_1(C^\infty(T^2)))\theta_1\wedge\bar\theta_1\}}=\infty,
	\]
	which implies   $\tilde h^{1,1}=\infty$.
	
	To calculate $\Tilde{H}^{2,0}_{Dol}$ and $\tilde{H}^{0,2}_{Dol}$ is not hard. Since the requirement 
	$u=f\theta_1\wedge\theta_2\in\ker(\bar\mu)$ shows that $f\equiv0$, i.e. $\tilde h^{2,0}=0$. Similarly, one also has $\Tilde{h}^{0,2}=0$.  
	
	Now, we calculate $\tilde{H}^{2,1}_{Dol}$ and $\tilde{H}^{1,2}_{Dol}$.
	For any $w=a
	\theta_1\wedge\theta_2\wedge\bar\theta_1+b\theta_1\wedge\theta_2\wedge\bar\theta_2\in\ker(\bar\partial)$, one has $-\bar v_2a+\bar v_1b=0$. By the identity $\ker(\bar\mu)\cap \mathcal{A}^{2,0}=\{0\}$, one has 
	\[\Tilde{H}^{2,1}_{Dol}\supset 
	\{w=a(t,x)\theta_1\wedge\theta_2\wedge\bar\theta_1\mid a\mbox{ depends only on }(t,x)\}, \]
	i.e., $\Tilde{h}^{2,1}=\infty$.
	
	\noindent 
	Since $\Tilde{H}^{1,2}_{Dol}=\frac{\mathcal{A}^{1,2}}{im(\bar\partial)}$, one has $\Tilde{H}^{1,2}_{Dol}\cong \ker(\bar\partial^*)\cap\mathcal{A}^{1,2}=\ker(\Delta_{\bar\partial})\cap\mathcal{A}^{1,2}$ by the Hodge decomposition of $\mathcal{A}^{1,2}$ with respect to $\Delta_{\bar\partial}$. Hence, we get 
	\[\Tilde{h}^{1,2}=\dim(\ker(\Delta_{\bar\partial}\cap\mathcal{A}^{1,2}))=\dim(\ker(\partial)\cap\mathcal{A}^{0,1})=\Tilde{h}^{1,0}=1.\]
	Similarly, one also has $\Tilde{h}^{2,2}=1$.
	Summarizing the above results, we write 
	\[\dim_\mbc\Tilde{H}^{*,*}_{Dol}(KT^4,J)=\begin{array}{ccccc}
		& &1&& \\
		&\infty &&1&\\
		0&&\infty&&0\\
		&1&&1&\\
		&&1&&
	\end{array}.\]
	
\end{exa}

\section{Almost complex 4 manifolds}

Throughout this section, all almost complex manifolds are compact and $4$ dimensional, unless otherwise specified.   Recall that for any Hermitian metric $g$ on an almost complex $4$ manifold, one has the identities, c.f., \cite[Lemma 2.1.57]{Don86},
\[\Lambda^+=\Lambda^{2,0}\oplus\Lambda^{0,2}\oplus\mathbb C\langle\omega\rangle,~\Lambda^-=\ker(L)\cap\Lambda^{1,1},\]
where $\omega$ is the imaginary part of $g$ and $\Lambda^+$ and $\Lambda^-$ are the $1$ and $-1$-eigenspaces of the Hodge-$*$ operator  with respect to $g$ respectively. 

\subsection{First Betti number on almost complex $4$ manifolds}
In this subsection, we show some inequalities about   the dimensions of Dolbeault cohomology groups  on compact almost complex $4$ manifolds, which are analogous to the case of compact complex surfaces.

\begin{lemma}\label{lemma-barpartial}
Let $(X,J)$ be  a compact  almost complex $4$ manifold. Then, it holds that 
\[\ker(\bar\partial)\cap \mathcal{A}^{1,0}=\ker(d)\cap \mathcal{A}^{1,0}.\]
\end{lemma}

\begin{pf}
For any form $\alpha\in \mathcal{A}^{1,0}(X)$, it si clear that
\[d\alpha=\partial\alpha+\bar\partial\alpha+\bar\mu\alpha.\]
The inclusion $\ker(d)\cap \mathcal{A}^{1,0}\subset \ker(\bar\partial)\cap \mathcal{A}^{1,0}$ is obvious.
Assume that $\bar\partial\alpha=0$. 
The  two identities $$*|_{\mathcal{A}^{2,0}}=Id \mbox{ and }*|_{\mathcal{A}^{0,2}}=Id,$$ 
shows that
\[d\alpha\wedge \overline{d\alpha}=\partial\alpha\wedge\overline{*\partial\alpha}+
\bar\mu\alpha\wedge\overline{*\bar\mu\alpha}
=|\partial\alpha|^2+|\bar\mu\alpha|^2.
\]
Integral the above formula   over $X$ gives that $\ker(\bar\partial)\cap \mathcal{A}^{1,0}\subset \ker(d)\cap \mathcal{A}^{1,0}$, which finishes the proof.
\end{pf}

Consequently, the above lemma and \cite[Lemma 2.3]{CW20} imply the following corollary.

\begin{cor}\label{cor-ell-h(1,0)}
On a compact almost complex $4$ manifold, it holds that
\[\ker(\Delta_{\bar\partial})\cap\ker(\Delta_{\bar\mu})
\cap\mathcal{A}
^{1,0}=\ker(\Delta_{\bar\partial}+\Delta_\mu)\cap \mathcal{A}^{1,0}=\ker(\bar\partial)\cap \mathcal{A}^{1,0}=\ker(\Delta_{\bar \partial})\cap\mathcal{A}^{1,0}.\]
\end{cor}

Similar to compact complex surfaces, one has the following lemmas.

\begin{lemma}\label{lemma-h1}
It holds that   $$0\to\tilde  H^{1,0}_{Dol}=H^{1,0}_{Dol}\to\tilde  H^{0,1}_{Dol}\subseteq\hat{H}^{0,1}\subseteq H^{0,1}_{Dol},$$
on any compact almost complex $4$ manifold.
\end{lemma}

\begin{pf}
The inclusions $\tilde  H^{0,1}_{Dol}\subseteq \hat H^{0,1}\subseteq H^{0,1}_{Dol}$ follow from the definitions. We just show the first map is an injection.
If an  element $u\in\ker(\bar\partial)\cap \mathcal{A}^{1,0}$  can be written as $u=\partial f$ for some  function $f\in\mathcal{A}^0_{Dol}$, then   $\bar\partial \partial f=0$. By applying the maximum principle to the equation $i\partial\bar\partial f=0 $,    $f$ is a constant function and $u=0$.  Hence, the first map is injective by the conjugation.
\end{pf}

\begin{lemma}\label{lemma-b1-equiality}
Let $(X,J)$ be a compact almost complex $4$ manifold. Then, it holds that $b_1\geq 2\tilde h^{1,0}=2h^{1,0}$. 
\end{lemma}

\begin{pf}
It is clear that $\tilde H^{1,0}_{Dol}\cap \overline{\tilde H^{1,0}_{Dol}}=\{0\}$ and $\tilde H^{1,0}_{Dol}=\ker(d)\cap\mathcal{A}^{1,0}$. The  fact  $$\tilde H^{1,0}_{Dol}\cap im(\partial)=\tilde H^{1,0}_{Dol}\cap im(d)=\{0\},$$ 
implies  the    desired inequality.
\end{pf}


\noindent
Now, we consider the relations on the Dolbeault and the refined Dolbeault  cohomologies with bi-degrees $(2,0)$ and $(0,2)$.
First, we give the following lemma.

\begin{lemma}\label{lemma-(2,0)-trivial}
Let $(X,J)$ be  a compact almost complex $4$ manifold.
If  $\sigma\in H^{2,0}_{Dol}$ can be represented by $\sigma=\partial \alpha+\mu\beta$ for some $\alpha\in \mathcal{A}^{1,0}$ and $\beta\in\mathcal{A}^{0,1}$, then $\sigma=0$.
\end{lemma}

\begin{pf}
Clearly,  the formulas $$\bar\partial\alpha\wedge\bar\sigma=0,~\bar\mu\alpha\wedge\bar\sigma=0,$$ and
\[\bar\partial\beta\wedge\bar\sigma=0, ~\partial\beta\wedge\bar\sigma=0,\] 
imply that
$$\sigma\wedge\bar\sigma=d(\alpha+\beta)\wedge\bar\sigma.$$ By   the Stokes lemma and the formula $d\sigma=\bar\partial\sigma+\bar\mu\sigma=0$, one obtains
\[
\int_X|\sigma|^2=\int_X\sigma\wedge \bar \sigma=0,
\]
i.e.,  $\sigma=0$.
\end{pf}

\noindent
The following lemmas  are  similar to the case of compact complex surfaces.

\begin{lemma}\label{lemma-iso-(2,0)-(0,2)}
The conjugation map $$H^{2,0}_{Dol}\to H^{0,2}_{Dol},$$
is an  isomorphism on a   compact almost complex $4$ manifold $(X,J,g,\omega)$.
\end{lemma}

\begin{pf}
Recall	the two isomorphisms in Theorem \ref{thm-CW-21} $$\ker(\Delta_{\bar\partial}+\Delta_{\bar\mu})\cap\mathcal{A}^{2,0}=H^{2,0}_{Dol}\mbox{ and }
\ker(\Delta_{\bar\partial}+\Delta_{\bar\mu})\cap\mathcal{A}^{0,2}\cong H^{0,2}_{Dol}.$$ 
By using 
the facts that the operators $\bar\partial$ and $\bar\mu$ act trivially on $\mathcal{A}^{0,2}$,  
$\sigma\mbox{ belongs to }\ker(\Delta_{\bar{\partial}}+\Delta_{\bar \mu})\cap\mathcal{A}^{0,2}$
if and only if $\sigma\in\ker(\bar\partial\bar\partial^*+\bar\mu\bar\mu^*)\cap\mathcal{A}^{0,2}$. 
For any $\sigma\in\ker(\Delta_{\bar{\partial}}+\Delta_{\bar \mu})\cap\mathcal{A}^{0,2}$, we get
\begin{eqnarray*}
	0=\int_X(\bar\partial\bar\partial^*\sigma+\bar\mu\bar\mu^*\sigma)\wedge\bar\sigma
	=\int_X|\bar\partial^*\sigma|^2+|\bar\mu^*\sigma|^2
	=\int_X|\partial\sigma|^2+|\mu\sigma|^2,
\end{eqnarray*}
where we use the self-duality of $\mathcal{A}^{0,2}$ for the last equality. On the other hand, the inclusion $$\ker(\bar\partial^*)\cap\ker(\bar\mu^*)\cap\mathcal{A}^{0,2}=\ker(\partial)\cap\ker(\mu)\cap\mathcal{A}^{0,2}\subseteq\ker(\Delta_{\bar{\partial}}+\Delta_{\bar \mu})\cap\mathcal{A}^{0,2},$$ 
is clear. 
Hence, by $\ker(\Delta_{\bar{\partial}})\cap\mathcal{A}^{0,2}=\ker(\partial)\cap\mathcal{A}^{0,2}$, 
$\ker(\Delta_{\mu})\cap\mathcal{A}^{0,2}=\ker(\mu)\cap\mathcal{A}^{0,2}$ and \cite[Lemma 2.3]{CW20}, the space of $(\Delta_{\bar{\partial}}+\Delta_{\bar \mu})$-harmonic $(0,2)$-forms is equal to the space $\mathcal{H}^{0,2}_{\bar\partial,\mu}$. 

By using the above arguments  to $(2,0)$-forms, 
it also holds 
$$H^{2,0}_{Dol}=\ker(\Delta_{\bar\partial}+\Delta_{\bar\mu})\cap\mathcal{A}^{2,0}=\ker(\bar\partial)\cap\ker(\bar\mu)\cap\mathcal{A}^{2,0}.$$
Thus, the conjugation
map
\[H^{2,0}_{Dol}\cong\ker( {\bar\partial})\cap \ker( {\bar\mu})\cap\mathcal{A}^{2,0}\to
\ker( {\partial})\cap \ker( {\mu})\cap\mathcal{A}^{0,2}
\cong H^{0,2}_{Dol},\]
is an isomorphism. This proves the lemma.
\end{pf}

\begin{lemma}
It holds that 
\[H^{2,0}_{Dol}=\tilde H^{2,0}_{Dol}=\mathcal{H}^{2,0}_{\bar\partial,\mu}
\cong\mathcal{H}^{0,2}_{\bar\partial,\mu}=H^{0,2}_{Dol}\subseteq \tilde{H}^{0,2}_{Dol},\]
on any compact almost Hermitian $4$ manifold.
\end{lemma}
\begin{pf}
By the first paragraph of the proof to Lemma \ref{lemma-iso-(2,0)-(0,2)}, one obtains 
$H^{0,2}_{Dol}\cong \ker(\partial)\cap\ker(\mu)\cap\mathcal{A}^{0,2}=\mathcal{H}^{0,2}_{\bar\partial,\mu}$. Similarly, we also get $\mathcal{H}^{2,0}_{\bar\partial,\mu}=\ker(\bar\partial)\cap\ker(\bar\mu)\cap\mathcal{A}^{2,0}$, which implies the isomorphism 
$\mathcal{H}^{2,0}_{\bar\partial,\mu}
\cong\mathcal{H}^{0,2}_{\bar\partial,\mu}$. This also gives a well-defined map 
$$H^{0,2}_{Dol}\cong\mathcal{H}^{0,2}_{\bar{\partial},\mu}\to \tilde{H}^{0,2}_{Dol}.$$ 
For any $\bar{\partial},\mu$-harmonic $(0,2)$-form $\sigma$, if it is trivial in $\tilde{H}^{0,2}_{Dol}$, then  $\sigma=\bar\partial\alpha$ for some $(0,1)$-form $\alpha$.  By taking the conjugation on $\sigma$ and using Lemma \ref{lemma-(2,0)-trivial}, we get $\sigma=0$, i.e., $\mathcal{H}^{0,2}_{\bar{\partial},\mu}\to \tilde{H}^{0,2}_{Dol}$ is an injection. 
\end{pf}

\noindent
We summarize the results of the above lemmas by the following corollary.

\begin{cor}\label{coro-betti-inequality}
Let $(X,J,g,\omega)$ be  a compact almost Hermitian $4$ manifold.
The following inequalities establish
\begin{itemize}
	\item[$(1)$] $\ell^{1,0}= h^{1,0}=\tilde h^{1,0}\leq \tilde h^{0,1} \leq \hat h^{0,1}\leq h^{0,1}$ and $ \ell^{0,1}  \leq \tilde h^{0,1} \leq\hat h^{0,1}\leq h^{0,1}$.
	\item[$(2)$] $\ell^{2,0}=h^{2,0}=\tilde{h}^{2,0}=\ell^{0,2}=
	h^{0,2}\leq\tilde{h}^{0,2}$.
\end{itemize}
\end{cor}
The first inequalities are the consequences of Lemma \ref{lemma-h1}, which prove \eqref{eqn-tilde-Dolbeault-ineq}. Here we only give the proof for $\ell^{0,1}\leq\tilde h^{0,1}$. 

\begin{pf}
We define  
$$\mathcal{H}^{0,1}_{\bar\partial}\cap \mathcal{H}^{0,1}_{\mu}\to\tilde H^{0,1}_{Dol},$$ 
by $u\mapsto [u]$. 
When $u$ maps to a trivial class, we write $u=\bar\partial f$ for some function $f\in\mathcal{A}^0_{Dol}$. Applying the maximum principle to the  equation $\bar\partial^*\bar\partial f=0$,   $f$ is a constant function and $u=0$.  
\end{pf}

\noindent
Now we show the relation between $b_1$ and $\tilde h^{1,0}+\hat h^{0,1}$. 

\begin{thm}\label{thm-b1-inequality}
	Let $(X,J)$ be a compact almost complex $4$ manifold. Then, it holds that $b_1\leq \tilde   h^{1,0}+\hat h^{0,1}\leq \tilde  h^{1,0}+  h^{0,1}$. 
	Moreover, when the equality holds $b_1=\tilde h^{1,0}+\hat h^{0,1}$, 
	we have that $\hat h^{0,1}=\dim_{\mathbb C}E^{0,1}_2$, where $E^{0,1}_2$ is the second stage of the Hodge filtration. 
\end{thm}


\begin{pf} 
	By the relation $\hat H^{0,1}\subseteq H^{0,1}_{Dol}$, we only need to prove the first inequality
$
b_1\leq\tilde  h^{1,0}+\hat h^{0,1}. 
$
Similar to the first paragraph in the proof of Lemma \ref{lemma-b1-h1}, we define a linear map 
\begin{equation}
	\ker(d)\cap\mathcal{A}^1 \to \mathcal{A}^{0,1},\label{eqn-b1-map}
\end{equation} 
$ \mbox{by } u'+u''\mapsto u''$, where $u'\in\mathcal{A}^{1,0}$ and $u''\in\mathcal{A}^{0,1}$. 

One can check that the above map  \eqref{eqn-b1-map} induces a   well-defined map  $$p:H^1_{dR}\to \hat H^{0,1}.$$ 
If an image $[u'']$ is trivial in $\hat H^{0,1}$, then we write $u''=\bar\partial f$ for some function $f$. Combining the equation
\[0=d(u'+u'')=\partial u'+\mu u''+\bar\partial u''+\bar\mu u'+\bar\partial u'+\partial u'',\] 
for some $u'\in\mathcal{A}^{1,0}$ with the two identities $\bar\partial^2+\bar\mu\partial=0\mbox{ and }\partial\bar\partial+\bar\partial\partial=0$ on functions, we obtain the formulas
\[\partial(u'-\partial f)=0,~\bar\mu(u'-\partial f)=0\mbox{ and }\bar\partial (u'-\partial f)=0\]
i.e., $u'-\partial f\in\tilde  H^{1,0}_{Dol}=H^{1,0}_{Dol}$. On the other hand, we have   a linear injection 
\[0\to \tilde{H}^{1,0}_{Dol}=H^{1,0}_{Dol}\to H^1_{dR},\]
defined by $v\mapsto [v+0]$.
Hence, we establish the short exact sequence $$0\to\tilde  H^{1,0}_{Dol}=H^{1,0}_{Dol}\to H^1_{dR}\to \hat H^{0,1},$$
via the identification $\ker(p)\cong\tilde  H^{1,0}_{Dol}$, which yields the inequality.
%

By the same arguments of the above paragraph, one also obtains an exact sequence
$$0\to\tilde  H^{1,0}_{Dol}= H^{1,0}_{Dol}\to H^1_{dR}\overset{p'}{\to }H^{0,1}_{Dol}.$$
When the equality $b_1=\tilde h^{1,0}+h^{0,1}$ holds, one gets   the exact sequence 
\begin{equation}
	0\to\tilde  H^{1,0}_{Dol}=H^{1,0}_{Dol}\to H^1_{dR}\to H^{0,1}_{Dol}\to0.
	\label{eqn-exact-b1}
\end{equation}
Recall that as in \cite[Appendix]{CW21} $E^{0,1}_2=E^{0,1}_\infty$ and $E^{1,0}_1=E^{1,0}_\infty$ for dimension $4$. By the equality $\dim E^{1,0}_\infty+\dim E^{0,1}_\infty=b_1$ of Theorem \ref{thm-CW-21}, one has the desired formula.
\end{pf}

\begin{rmk}
When $(X,J)$ is a compact  complex surface,  one can establish a similar inequality   by using the sheaf theory,  c.f., \cite[Lemma 2.4-Chapter IV]{BHPV}. 
\end{rmk}

By Theorem \ref{thm-b1-inequality}, we  have the following corollary.
\begin{cor}\label{cor-b1-ineq-4-mfld}
The inequalities 
\[2\tilde h^{1,0}=2h^{1,0}\leq b_1\leq \tilde h^{1,0}+\hat h^{0,1}\leq \tilde h^{1,0}+h^{0,1},\]
hold on any compact almost complex $4$ manifold. Moreover, one has the following: 
\begin{itemize}
	\item[$(1)$] The equality $\tilde h^{1,0}=\hat{h}^{0,1}$ implies 
	$
	b_1=\hat h^1=2\tilde h^{1,0}
	$.
	\item[$(2)$]  The equality   $b_1=\hat h^1$ implies 
	$\tilde h^{1,0}=\tilde h^{0,1}.$
	\item[$(3)$] The equality  $\tilde h^{1,0}=h^{0,1}$ implies the identities  
	$h^{1,0}=\ell^{1,0}=\hat h^{0,1}=\tilde{h}^{0,1}.$
\end{itemize}

\end{cor}

\begin{pf}
The three inequalities follows from Theorem \ref{thm-b1-inequality} and Lemma \ref{lemma-b1-equiality}. 
\begin{itemize}
	\item[$(1)$]  The equality $\tilde h^{1,0}=\hat h^{0,1}$ is equivalent to the  condition $\tilde h^{1,0}=\tilde h^{0,1}=\hat{h}^{0,1}$.  Thus,  the result follows from the inequalities
	$$2\tilde h^{1,0}\leq b_1\leq \tilde h^{1,0}+\hat h^{0,1}\leq \tilde h^{0,1}+\hat h^{0,1}=\hat h^1.$$
	\item[$(2)  $]   Recall that there are  two inequalities
	$$b_1\leq\tilde h^{1,0}+\hat{h}^{0,1}\mbox{ and }\tilde h^{1,0}\leq \tilde{h}^{0,1}.$$ 
	Hence, the equality $b_1= \hat h^1=\tilde{h}^{0,1}+\hat{h}^{0,1}$ shows $	\tilde h^{1,0}=\tilde{h}^{0,1}$. 
	\item[$(3)$] It follows from Corollary \ref{coro-betti-inequality}.
\end{itemize}
\end{pf}


On compact almost K\"ahler manifolds, Cirici and Wilson \cite{CW20} gave  the equality $\ell^{p,q}=\ell^{q,p}$ by using  the relation $$\Delta_{\bar\partial}+\Delta_{\mu}=\Delta_\partial+\Delta_{\bar\mu},$$ and the  complex conjugate relation in Theorem \ref{thm-CW-02}. 
For the completeness of this paper, we  show the equality $\ell^{1,0}=\ell^{0,1}$ on compact almost K\"ahler $4$ manifolds via a simple method rather than using $\bar\partial,\mu$-Laplacian operator. 
\begin{lemma}
On a compact almost K\"ahler $4$ manifold $(X,J,g,\omega)$, 
it holds that
$\tilde h^{1,0}=\ell^{1,0}=\ell^{0,1}$.
\end{lemma}
\begin{pf} 
The first equality is provided by Corollary \ref{cor-ell-h(1,0)}. We just show the second one here. The formula 	$[\Lambda, \bar{\partial}]=-i \partial^*$ in Proposition \ref{prop-kahler-identities} implies $\partial^*v=0$ for $v\in\tilde{H}^{1,0}_{Dol}$, which gives an injection $0\to\tilde{H}^{1,0}_{Dol}\to \mathcal{H}^{0,1}_{\bar\partial,\mu}$ by the conjugation map.

Similarly, the formula $[\Lambda, \partial]=i \bar{\partial}^*$ in Proposition \ref{prop-kahler-identities} induces  that
$i\bar\partial^* u=\Lambda(\partial u)$ for any $u\in \mathcal{A}^{0,1}$. 
For any $u\in\mathcal{H}^{0,1}_{\bar\partial,\mu}$, the conditions $\bar\partial^*u=0$ and $\bar\partial u=\mu u=0$  imply that $d^+u=0$. Since on a compact oriented $4$ smooth manifold without boundary,  the formulas
\[0=\int_Xdu\wedge d\bar u=\int_Xd^+ u\wedge d^+\bar u+\int_Xd^-u\wedge d^-\bar u=\int_X|d^+u|^2-|d^-u|^2,\]
and $\|du\|^2_{L^2}=\int_X|d^+u|^2+|d^-u|^2$ hold.  
Hence,  we get that  $du=0$ and  the conjugation map $\tilde{H}^{1,0}_{Dol}\to \mathcal{H}^{0,1}_{\bar\partial,\mu}$ is an isomorphism. 
\end{pf}

Recall that
on a compact symplectic $4$ manifold, its first Betti number is not even  by Gompf's theorem, in other words the equality $b_1=2\tilde h^{1,0}$ does not hold for a general compact  almost K\"ahler $4$ manifold. At the end of this subsection, we give an analytic condition for a compact  almost K\"ahler $4$ manifold admitting an even first Betti number.

\begin{lemma}\label{lemma-b1-almost-kahler}
Suppose that   the inequality
\[
\lambda_1>4\|\Delta_{\bar\mu}|_{\mathcal{A}^1}\|,\]
holds on a compact almost K\"ahler $4$ manifold $(X,J,g,\omega)$,
where $\lambda_1$ denotes the first non-zero eigenvalue of the Hodge-Laplacian operator $\Delta_d$ restricted on $1$-forms and 
$\|\Delta_{\bar\mu}|_{\mathcal{A}^1}\|$ denotes the operator norm of $\Delta_{\bar \mu}:L^2(\mathcal{A}^1)\to L^2(\mathcal{A}^1)$.
Then,  $b_1=2\tilde h^{1,0}=2\ell^{1,0}=2\ell^{0,1}$.
\end{lemma}
\begin{pf}
Recall that $\Delta_d$ is a real operator, and each real-valued $1$-form in $\mathcal{A}^1_{\mathbb R}$ can be written as $u+\bar u$ for a unique element $u\in\mathcal{A}^{1,0}$, where $\mathcal{A}^1_{\mathbb R}$ is the space of real-valued $1$-forms.  By Lemma \ref{lemma-b1-equiality}, it suffices to show that each element $u+\bar u\in\ker(d)\cap\ker(d^*)\cap\mathcal{A}^1_{\mathbb R}$ corresponds to a unique element $u\in\tilde H^{1,0}_{Dol}$ on $(X,J,g,\omega)$.  

Combining  the condition $d^*(u+\bar u)=0 ,$   
\mbox{with the K\"ahler identities  }$[\Lambda,\partial]=i\bar\partial^*\mbox{ and }[\Lambda,\bar\partial]=-i\partial^*$
in Proposition \ref{prop-kahler-identities}, one obtains  $\Lambda(\bar\partial u-\partial\bar u)=0$. 
Together with the equation  $(\bar\partial u+\partial \bar u)=0$, we get  
\begin{equation}
	\Lambda\bar\partial u=\Lambda\partial \bar u=0,\label{eqn-asd}
\end{equation}
i.e., $\partial\bar u$ is an imaginary-valued anti-self-dual $2$-form. 
We have the formula   
\begin{eqnarray*}       
	\int_X\partial u\wedge\bar\partial\bar u&=&        -\int_X\bar\partial\partial u\wedge\bar u\\        &=&-\int_X(\partial^2\bar u+\mu\bar\partial\bar u)\wedge \bar u\\       
	&=&\int_X\partial\bar u\wedge\partial\bar u+\int_X\bar\partial\bar u\wedge\mu\bar u\\        &=&\int_X\partial\bar u\wedge\partial\bar u-\int_X\bar\partial\bar u\wedge \partial u,   
\end{eqnarray*}   
where we use  the formula $\bar\partial\partial u=-\partial\bar\partial u-\bar\mu\mu u-\mu\bar\mu u=    \partial^2\bar u+\mu\bar\partial\bar u$ for the second equality and the formula $\partial u+\mu\bar u=0$ for the last one.   
The above formula and the anti-self-duality of $\bar\partial u$ yield 	\begin{equation}
	\|\bar\partial u\|^2_{L^2}=\int_X\partial\bar u\wedge*\bar\partial u=-\int_X\bar\partial u\wedge\partial\bar u=\int_X\partial\bar u\wedge\partial\bar u=2\|\partial u\|^2_{L^2}.\label{eqn-parital-d}
\end{equation}  
On the other hand, we compute  
\begin{eqnarray*}     
	\int_X(\Delta_du, u)&=&\int_X(du, du) +(d^*u,d^*u) \\
	&=&\int_X(\bar\partial u,\bar\partial u)+(\partial u,\partial u)+(\bar\mu u,\bar\mu u)  = 4\int_X(\Delta_{\bar\mu} u,u), 
\end{eqnarray*}    
where we use the formulas $ d^*u=\partial^*u=0$, \eqref{eqn-asd} for the second equality and  the formulas $\bar\partial \bar u+\bar\mu u=0$, \eqref{eqn-parital-d} for the last one. Therefore, by the hypothesis of $\lambda_1$, it holds that $u\in\ker(d)\cap\mathcal{A}^{1,0}$, i.e., $u\in\tilde  H^{1,0}_{Dol}$.
\end{pf}

\subsection{Symplectic condition on almost complex $4$ manifolds}

First, we show that the equality $\tilde h^{1,0}=\tilde h^{0,1}$ is equivalent to the condition that each $dd^c$-closed $(1,1)$-form generates a unique  $d$-closed $2$-form. 
For the $(2,0)$ and $(0,2)$-components of the form
$d(u+\bar u)=\partial u+\bar\partial \bar u+\mu u+\bar\mu \bar u+\bar\partial u+\partial \bar u$, one has  
$$d(\partial\bar u+\bar\partial u)=\bar\partial\partial \bar u+\partial\bar\partial  u+\bar\mu\partial\bar u+\mu\bar\partial u,$$
and
\[
d(\mu u+\bar\mu\bar u)=\bar\mu\mu u+\mu\bar\mu\bar u+\bar\partial \mu u+\partial\bar\mu\bar u
\]where $u\in\mathcal{A}^{0,1}$.
For a $dd^c$-closed real $(1,1)$-form $\psi$, we want to find some  $u\in\mathcal{A}^{0,1}$  satisfying the equation,
$$d\psi=d(\bar\partial u+\partial\bar u+\mu u+\bar\mu\bar u).$$ 
This is equivalent to solving $\bar\partial\psi=\partial\bar\partial u+\bar\mu\partial\bar u+\partial\bar\mu \bar u+\bar\mu\mu u$. The existence of such $u$   is provided by the following lemma.  


\begin{lemma}\label{lemma-partial-solution}
Let $(X,J,g,\omega)$ be a compact almost Hermitian $4$ 
manifold. Then, the condition $\tilde h^{1,0}=\tilde h^{0,1}$ holds  if and only if for any $\psi\in\mathcal{A}^{1,1}_{\mathbb R}$ satisfying $\partial\bar\partial\psi=0$, there exists a solution  $u\in\mathcal{A}^{0,1}$ to the equation
\begin{equation}
	\bar\partial \psi=\partial\bar\partial u+\bar\mu\partial\bar u+\partial\bar\mu \bar u+\bar\mu\mu u,\label{eqn-ddbar}
\end{equation}
where $\mathcal{A}^{1,1}_{\mathbb R}$ denotes the space of real-valued $(1,1)$-forms.
\end{lemma}
\begin{pf}
Consider the adjoint operator,
\begin{eqnarray*}
	(\partial\bar\partial +\bar\mu\mu+(\partial\bar\mu+\bar\mu\partial)\comp c)^*&=&*(\partial\bar\partial)*+*(\bar\mu\mu)*+c *(\bar\partial\mu+\mu\bar\partial)*\\
	&=&*((\partial\bar\partial)+(\bar\mu\mu)+c\comp  (\bar\partial\mu+\mu\bar\partial))*,
\end{eqnarray*}
where $c$ denotes the conjugation map.
The equation has a solution if and only if the form $\bar\partial\psi$ is vertical to $\ker((\partial\bar\partial)*+(\bar\mu\mu)*+c\comp(\bar\partial\mu+\mu\bar\partial)*)$. 
This is equivalent to the condition that the formula
$\int_X  *\bar w\wedge \bar\partial\psi=0$
holds  
for all $w\in\mathcal{A}^{1,2}$ with the property $\partial\bar\partial*w+{\bar\mu\mu *w}+c\comp(\bar\partial\mu+\mu\bar\partial)*w=0$.
Setting $v=*  w$, 
we rewrite
\begin{equation}
	\partial\bar\partial v+\bar\mu\mu v+c\comp(\mu\bar\partial+\bar\partial\mu)v=0.\label{eqn-ddbar-c}
\end{equation}
By the Stokes lemma, the above equation \eqref{eqn-ddbar-c} implies the following:
\begin{eqnarray*}
	0&=&\int_X (\partial\bar\partial v+\bar\mu\mu v+c\comp(\mu\bar\partial+\bar\partial\mu)v)\wedge \bar v\\
	&=&
	\int_Xd(\bar\partial v+\mu v+\partial \bar v+\bar\mu\bar v)\wedge\bar v\\
	&=&-\int_X(\bar\partial v+\mu v+\partial \bar v+\bar\mu\bar v)\wedge d(\bar v)\\
	&=&-\int_X (\bar\partial v+\mu v+\partial \bar v+\bar\mu\bar v)\wedge 
	(\partial\bar v+\bar\mu\bar v)\\
	&=&\int_X-(\bar\partial v,\bar\partial v)-(\bar\mu \bar v,\bar\mu \bar v)-(\bar\partial  v,\bar\mu \bar v)-(\bar\mu \bar v,\bar\partial v)\\
	&=& 
	-\|\bar\partial v+\bar\mu\bar v\|^2_{L^2},
\end{eqnarray*}
i.e.,  $[v]\in\hat{H}^{0,1}$ and $[v+\bar v]\in \hat{H}^1$.
The condition $\tilde h^{1,0}=\tilde h^{0,1}$ establishes an  exact sequence $$0\to\tilde  H^{1,0}_{\partial}\cong\tilde{H}^{1,0}_{Dol}\to\hat H^1\to\hat H^{0,1}\to0.$$
Recall that $d^{2,0}$ and $d^{0,2}$ denote the $(2,0)$ and $(0,2)$ components of $im(d)$ respectively. 
By the above exact sequence and the condition $v+\bar v\in \ker(d^{0,2})\cap\ker(d^{2,0})$, there are two functions $f$ and $f'$ with the property $$d^{2,0}(\partial f+\bar\partial f')=\partial^2(f'-f'')=0,~d^{0,2}(\partial f+\bar\partial f')=\bar\partial^2(f''-f')=0,$$
such that
the $(1,0)$-component of $v+\bar v-\partial f-\bar\partial f'$ belongs to $\tilde{H}^{1,0}_{Dol}=\ker(d)\cap\mathcal{A}^{1,0}$, i.e., $\bar v-\partial  f\in\ker(d)$.
The equation $dd^c\psi=i\partial\bar\partial\psi=0$  implies that
\[\int_X\partial f\wedge\bar\partial\psi=-\int_Xf\wedge\partial\bar\partial\psi=0.\]
Therefore, it holds that
\[\int_X \bar v \wedge \bar\partial \psi
=\int_Xd(\bar v-\partial f) \wedge \psi=0,\]
i.e., $\int_X *\bar w\wedge \bar\partial\psi=0$. 

Conversely,  assume that for any $dd^c$-closed form $\psi\in\mathcal{A}^{1,1}_{\mathbb R}$  the equation $$\bar\partial \psi=\partial\bar\partial u+\bar\mu\partial\bar u+\bar\mu\partial \bar u+\bar\mu\mu u,$$
can be solved for some $u\in\mathcal{A}^{0,1}$. Recall that for any $\partial\bar\partial$-closed $(1,1)$-form $\psi$, it can be written as $\psi=\psi'+i\psi''$, where $\psi'$ and $i\psi''$ are the real  and the imaginary parts of $\psi$ respectively.
Since the operator $i\partial\bar\partial$ is real, the $\partial\bar\partial$-closeness of $\psi$ is equivalent to  those  of both $\psi'$ and $\psi''$.
Take the sum of the solutions $u'$ and $u''$ to the equations $$	\bar\partial \psi'=\partial\bar\partial u'+\bar\mu\partial\bar u'+\partial\bar\mu \bar u'+\bar\mu\mu u'\mbox{ and }	\bar\partial \psi''=\partial\bar\partial u''+\bar\mu\partial\bar u''+\partial\bar\mu \bar u''+\bar\mu\mu u''\mbox{ respectively}.$$ 
The assumption implies that for any $dd^c$-closed  $(1,1)$-form $\psi$, there are  solutions $u_1\in\mathcal{A}^{0,1}$ and $u_2\in\mathcal{A}^{1,0}$ to the equation 
\[\bar\partial\psi=\partial\bar\partial u_1+\bar\mu\mu u_1+\bar\mu\partial u_2+\partial\bar\mu u_2.\]
We need to show that each solution to the equations $\partial v=0$ and $\bar\mu v=0$ for $v\in\mathcal{A}^{1,0}_{Dol}$ generates a unique form in $\ker(\bar\partial)\cap\mathcal{A}^{1,0}$. 
Similar to the first paragraph, 
the equation $$\bar\partial v=\bar\partial\partial f,$$
can be solved for some function $f$, if and only if the formula $\int_X \bar\partial v\wedge \psi=0$ holds for any $dd^c$-closed $(1,1)$-form $\psi$. 
Again by the Stokes lemma,   the above formula is rewritten as follows:
\begin{eqnarray*}
	\int_X\bar\partial v\wedge \psi&=&\int_X v\wedge\bar\partial \psi\\
	&=&\int_X v\wedge(\partial\bar\partial u_1+\partial\bar\mu u_2+\bar\mu\partial u_2+\bar\mu\mu u_1)\\
	&=&\int_X
	(\partial v\wedge ({\bar\partial u_1+\bar\mu u_2})+\bar\mu v\wedge(\partial u_2+\mu u_1))=0.
\end{eqnarray*}
Consequently,   by the conjugation map,  each element in $\ker(\bar\partial)\cap\ker(\mu)\cap\mathcal{A}^{0,1}_{Dol}$  is $\partial$-closed up to some element in $im(\bar\partial)\cap\mathcal{A}^{0,1}_{Dol}$, i.e.,  $\tilde h^{1,0}=\tilde{h}^{0,1}$. 
\end{pf}

\begin{rmk}
For the case of compact complex surfaces, i.e. $\mu=\bar\mu=0$, Buchdahl \cite[Lemma 8]{Buch} considered the    equation for complex-valued $(1,1)$-form $\psi$.
\end{rmk}
Before proceeding, we verify the following  two relations on almost complex $4$ manifolds.
\begin{lemma}\label{lemma-ddbar-zero}
Let  $(X,J)$ be an almost complex $4$ manifold. We have the following formulas: 
\begin{itemize}
\item[$(1)$] $\partial\bar\partial (\partial u+\bar\partial v)=0$, for any $u\in\mathcal{A}^{0,1}$ and $v\in\mathcal{A}^{1,0}$;
\item[$(2)$] $\partial\bar\partial(\partial \bar\partial f)=0$, for any function $f$. 
\end{itemize}
\end{lemma}
\begin{pf}
From the formula \[\partial\bar\partial f=\partial(\bar\partial f/2)+\bar\partial(-\partial f/2),\] we see that (1) implies (2). Therefore, it suffices  to show (1). 
We formulate 
\begin{eqnarray*}\partial\bar\partial (\partial u+\bar\partial v)&=&-
\bar\partial\partial^2 u+\bar\partial\partial\bar\partial v=
(\bar\partial \bar\partial \mu u+\bar\partial \mu\partial u)+
(\bar\partial^2\partial v+\bar\partial \mu\bar\mu v)\\
&=&(\bar\partial^2 d(u+v)+\bar\partial \mu(d(u+v)))\\
&=&  -(\bar\partial^2+\bar\partial\mu+\mu\bar\partial)(d(u+v))=0.
\end{eqnarray*}
Here we use  the fact that $\mu,~\bar\mu $ act trivially on $\mathcal{A}^{1,1}$ for the first equality, $\bar\partial^2$ acts trivially on $\mathcal{A}^{1,1}$ and $\mathcal{A}^{0,2}$, $\bar\partial\mu$ acts trivially on  $\mathcal{A}^{1,1}$ and $\mathcal{A}^{2,0}$ for the second one, and $\mu$ acts trivially on $\mathcal{A}^3$ for the third one. 
\end{pf}

\noindent
Combining the above Lemma \ref{lemma-ddbar-zero} with Lemma \ref{lemma-partial-solution}, one has the following. 
\begin{lemma}\label{lemma-s-map}
Let $(X,J,g,\omega)$ be a compact almost Hermitian $4$ manifold.     
The condition  $\tilde h^{1,0}=\tilde h^{0,1}$ holds, if and only if there exists a linear map $s:Z^{1,1}_{dd^c}(X,\mathbb R)\to Z^{2}_{d}(X,\mathbb R)$, such that
\begin{enumerate}
\item $s(\psi)=\psi+\bar\partial u+\partial\bar u+\mu u+\bar\mu\bar u$, where $u\in\mathcal{A}^{0,1}$;
\item $s(\partial u+\bar\partial \bar u)=d(u+\bar u)$, for any $u\in \mathcal{A}^{0,1}$;
\end{enumerate}
where $Z^{1,1}_{dd^c}(X,\mathbb R)$ and $Z^2_d(X,\mathbb R)$ are the spaces of all $dd^c$-closed real-valued $(1,1)$-forms and    $d$-closed real-valued $2$-forms respectively. 
\end{lemma}
\begin{pf}
It suffices to prove that  any two solutions to the equation  \eqref{eqn-ddbar} induce the same $d$-closed $2$-form.

Assume that there are two forms $u_1,u_2\in\mathcal{A}^{0,1}$ satisfying  
\[\bar\partial \psi=\partial\bar\partial u_i+\bar\mu\partial\bar u_i+\partial\bar\mu \bar u_i+\bar\mu\mu u_i,\]
for some $dd^c$-closed real-valued $(1,1)$-form $\psi$ with $i=1,2$. 
Setting $\gamma=u_1-u_2$, it holds that
\begin{equation}\partial\bar\partial \gamma+\bar\mu\partial\bar \gamma+\partial\bar\mu \bar \gamma+\bar\mu\mu \gamma=0.\label{eqn-ddbar-zero}
\end{equation}
By the elementary identity $Re(a+b)=Re(a+\bar b)$ for any complex numbers $a$ and $b$, 
we compute
\begin{eqnarray*}
\int_X|\partial\bar\gamma+\mu\gamma|^2
&=&
Re\int_X(|\partial\bar\gamma|^2+|\mu\gamma|^2)+Re(\int_X(\partial\bar\gamma,\mu\gamma)+(\mu\gamma,\partial\bar\gamma))\\
&=&   Re\int_X\partial\bar\gamma\wedge\bar\partial \gamma+\mu\gamma\wedge\bar\mu\bar\gamma+
Re \int_X(\partial\bar\gamma\wedge\bar\mu\bar\gamma+
\overline{(\mu\gamma\wedge\bar\partial\gamma)}\\
&=&-Re(\int_X\partial\bar\partial\gamma\wedge\bar\gamma+
\bar\mu\mu\gamma\wedge\bar\gamma+ \bar\mu\partial\bar\gamma\wedge\bar\gamma+\partial\bar\mu\bar\gamma\wedge\bar\gamma)=0.
\end{eqnarray*}
This shows $\bar\partial\gamma+\bar\mu\bar\gamma=0=\partial\bar\gamma+\mu\gamma$. Substituting $\gamma=u_1-u_2$ into these  two equations, it is clear that \[
\bar\partial u_1+\bar\mu \bar u_1=\bar\partial u_2+\bar \mu\bar u_2
\mbox{ and }
\partial\bar u_1+\mu u_1=\partial\bar u_2+\mu u_2.
\]
Therefore, it yields $$
\bar\partial u_1+\partial\bar u_1+\mu u_1+\bar\mu\bar u_1=\bar\partial u_2+\partial\bar u_2+\mu u_2+\bar\mu\bar u_2,
$$ which shows that the map $s$ is well-defined. 
\end{pf}

\noindent
By Lemma \ref{lemma-ddbar-zero} and Lemma \ref{lemma-s-map},
one    has the following corollary.
\begin{cor}
Let $(X,J)$ be a compact almost complex $4$ manifold. Suppose that $\tilde h^{1,0}=\tilde h^{0,1}$. Then, the above map $s$ descends to an injective map
\[
S: H^{1,1}_{dd^c}(X;\mathbb R)\to H^2_{dR}(X;\mathbb R),
\]
where $H^{1,1}_{dd^c}:=\frac{\ker(dd^c:\mathcal{A}^{1,1}_{\mathbb R}\to\mathcal{A}^4_{\mathbb R})}{im(d^{1,1}:\mathcal{A}^1_{\mathbb R}\to\mathcal{A}^{1,1}_{\mathbb R})}$.
\end{cor}

Now we are ready to show our main theorem.

\begin{pf}\textbf{ of Theorem \ref{thm-main-1}}
By Gauduchon's result  \cite{Gau77}, one can find an Hermitian metric $g$ on $(X,J)$, whose imaginary part $\omega$  is  $dd^c$-closed. Lemma \ref{lemma-partial-solution} provides that there is a form $u\in\mathcal{A}^{0,1}$ such that the real-valued $2$-form $\omega'=\omega+\bar\partial u+\partial\bar u+ \mu u+\bar\mu\bar u$ is $d$-closed. 
By using the following Lemma \ref{lemma-nondegenerate},    $\omega'$ is $J$-taming, i.e., $(X,\omega')$ is a compact $J$-taming symplectic $4$ manifold. 
\end{pf}

\begin{lemma}\label{lemma-nondegenerate}
Let $\psi$ be a real-valued positive $(1,1)$-form of an almost complex $4$ manifold $(X,J)$ and let $\sigma$ be a $(2,0)$-form. Then, the real form 
$\psi'=\psi+\sigma+\bar\sigma$ is non-degenerate. 
\end{lemma}
\begin{pf}
We consider the Hermitian metric $g=\psi(J,)$. For a point $x$ of $X$, we choose a unitary basis $(Z_1,Z_2)$  with respect to $g$. Locally, one has the expressions 
\[\psi=i\sum_{k=1,2}\theta^k\wedge\bar\theta^k\mbox{ and }\sigma=a\theta^1\wedge\theta^2,\] 
where $(\theta^1,\theta^2)$ is dual to $(Z_1,Z_2)$ at $x$.
It suffices to show that for any  real-valued vector field $v\in\Gamma(TX)$ the contraction $\iota_v\psi'$ is zero at $x$ if and only if $v$ is zero at $x$. 

Again, we have the local expressions \[v=v_1Z_1+v_2Z_2+\bar v_1\bar Z_1+\bar v_2\bar Z_2,\]
and 
$$\iota_v\psi'=av_1\theta^2+iv_1\bar\theta^1-
av_2\theta^1+iv_2\bar\theta^2-i\bar v_1\theta^1+\bar a\bar v_1\bar\theta^2-i\bar v_2\theta^2-\bar a\bar v_2\bar\theta^1.$$
Therefore, the form $\iota_v\psi'$ is zero at $x$  if and only if
\[
av_1-i\bar v_2=0,~i  v_1-\bar a\bar v_2=0.
\]
The non-singularity of the matrix $\left(\begin{array}{cc}
a&-i  \\
i & -\bar a
\end{array}\right)$ shows that $v_1=\bar v_1=v_2=\bar v_2=0$, i.e., $v$ is zero at $x$. 
\end{pf}

\hfill

Now, we show the generalized  $\partial\bar\partial$-lemma on compact almost complex $4$ manifolds without assuming the vanishing of Nijenhuis tensor. 

\begin{pf} {\bf of Theorem \ref{thm-ddbar}}
First,   assume that the equality $\tilde h^{1,0}=\tilde h^{0,1}$ holds. Recall that this establishes the exact sequence
\[0\to \tilde{H}^{1,0}_{\partial}\cong\tilde H^{1,0}_{Dol}\to \hat H^1\to \hat H^{0,1}\to0.	\] 
For any $d$-exact $(1,1)$-form $\psi$,   we write $\psi=d\alpha$ for some $\alpha\in \mathcal{A}^1$.  
The condition that $d\alpha$ is a $(1,1)$-form shows that the $(2,0)$ and $(0,2)$ components of $d\alpha$ vanish, i.e.,
 \[\bar\partial \alpha''+\bar\mu\alpha'=0\mbox{ and }\partial\alpha'+\mu\alpha''=0,\]
where $\alpha'$  and $\alpha''$ are  the $(1,0)$ and  $(0,1)$-component of $\alpha$ respectively. 

Since
$[\alpha'+\alpha'']\in\hat H^{1}$, by definition there are two functions  $ f'$ and $f''$ satisfying 
\begin{equation}
 (\partial^2 f'-\partial^2 f'')=0= (\bar\partial^2 f''-\bar\partial^2 f'),\label{eqn-d-0}
\end{equation}
such that
\begin{equation}
\alpha'-\partial f' \in\ker(d)\cap\mathcal{A}^{1,0}=\tilde{H}^{1,0}_{Dol}.\label{eqn-d-1}
\end{equation}
Together with the equations \[\bar\partial (\alpha''-\bar\partial f'')+\bar\mu (\alpha'-\partial f')=0,\]
and
\[\partial (\alpha'-\partial f')+\mu (\alpha''-\bar\partial f'')=0,\] it derives 
\[\alpha''-\bar\partial f''\in\ker(\mu)\cap\ker(\bar\partial),\]
i.e., $[\alpha''-\bar\partial f'']\in\tilde H^{0,1}_{Dol}$.

Again,   the condition $\tilde h^{1,0}=\tilde{h}^{0,1}$ implies that there is a function 
\begin{equation}
f\in\ker(\partial^2)\cap\ker(\bar\partial^2),\label{eqn-d-f}
\end{equation}
such that 
\begin{equation}
\alpha''-\bar\partial f''-\bar\partial f\in\ker(d)\cap\mathcal{A}^{0,1}. \label{eqn-d-2}
\end{equation}
Hence, we obtain
\begin{eqnarray*}
\psi&=&d(\alpha''-\bar\partial f''-\bar\partial f)+d(\alpha'-\partial f')+d(\bar\partial f''+\partial f'+\bar\partial f)\\
&=&d^{2,0}(\bar\partial f''+\partial f'+\bar\partial f)+d^{0,2}(\bar\partial f''+\partial f'+\bar\partial f)+d^{1,1}(\bar\partial f''+\partial f'+\bar\partial f)\\
&=&\partial \bar\partial(f''-f'+f) ,
\end{eqnarray*}
where we use \eqref{eqn-d-1}, \eqref{eqn-d-2} for the second equality and \eqref{eqn-d-0}, \eqref{eqn-d-f} for the last equality.   


Conversely,    assume that each $d$-exact  $(1,1)$-form $\psi $ can be  written as $\psi=\partial\bar\partial f$ for some   function $f$. We consider a representative form $u$ in a class $[u]\in \tilde H^{0,1}_{Dol}$. Clearly, the form $\partial u=du$ is a $d$-exact $(1,1)$-form. By the assumption, there is a function $f$ such that 
$$\partial (u-\bar\partial f)=0.$$
We obtain $u-\bar\partial f\in \ker(d)\cap\mathcal{A}^{0,1}$ by Lemma \ref{lemma-barpartial}, which implies that $f\in\ker(\partial^2)\cap\ker(\bar\partial^2)$ by $\mu(u)=\bar\partial(u)=0$. 
Therefore,  the conjugation map 
\[ 0\to  H^{1,0}_{Dol}=\tilde H^{1,0}_{Dol}\to \tilde H^{0,1}_{Dol},\]
is an isomorphism, i.e., $\tilde h^{1,0}=\tilde h^{0,1}$.  
\end{pf}


\noindent
One also has the following corollary. 
\begin{cor}
Let $(X,J)$ be a compact almost complex $4$ manifold. Then, the condition $\tilde h^{1,0}=\tilde h^{0,1}$ is equivalent to the 
condition
$$H^{1,1}_{dR}(X;\mathbb C)\cong H^{1,1}_{BC}(X;\mathbb C),$$ where 
$H^{1,1}_{BC}(X;\mathbb C):=\frac{\ker(d)\cap\mathcal{A}^{1,1}}{
im(dd^c)\cap\mathcal{A}^{1,1}
}$ and 
$H^{1,1}_{dR}(X;\mathbb C):=\frac{\ker(d)\cap\mathcal{A}^{1,1}}{
im(d)\cap\mathcal{A}^{1,1}
}$.
\end{cor}

At the end of this paper, we give an application. 
Recall that on compact almost complex $4$ manifolds, we show the following relations, c.f., Corollary \ref{cor-b1-ineq-4-mfld}:\begin{equation}	2\tilde{h}^{1,0}=2h^{1,0}\leq b_1 \mbox{ and }\tilde{h}^{1,0}\leq\tilde{h}^{0,1}\leq\hat{h}^{0,1}\leq h^{0,1}.	\label{eqn-b1-ineq-4mfld}\end{equation}  
Our definition of the  refined Dolbeault cohomology can also be used to detect some difference between the compact complex surfaces and the general almost complex $4$ manifolds.
It is clear   that on compact complex surfaces  the equality 
\begin{equation}
	b_1=h^{1,0}+h^{0,1},\label{eqn-b1-equality}
\end{equation}
always holds by using the Hodge decomposition and the Atiyah-Singer-Hirzebruch index theory, c.f., \cite[Theorem 2.7-Chapter IV]{BHPV}. By  definition,  the equalities $\tilde h^{0,1}=\hat{h}^{0,1}=h^{0,1}$ also trivially hold on  compact complex surfaces. Hence, it is natural to wonder  whether an analog of \eqref{eqn-b1-equality} in terms of $\tilde{h}^{1,0},~\tilde{h}^{0,1},~\hat{h}^{0,1},~h^{0,1}$ and $b_1$ holds for compact  almost complex $4$ manifolds. However, by Theorem \ref{thm-main-1} and Inequalities \eqref{eqn-b1-ineq-4mfld}, neither of  the equalities
\begin{equation}
	b_1=\tilde h^{1,0}+ \tilde h^{0,1},~b_1=\tilde h^{1,0}+ \hat h^{0,1},\mbox{ }b_1=\tilde  h^{1,0}+{h}^{0,1}, \label{eqn-b1-hodge-complex-surface}
\end{equation}  holds on general compact  almost complex $4$  manifolds. Moreover, this implies that the three identities in \eqref{eqn-three-identities} are not equivalent to each other on compact almost complex $4$ manifolds. 
For example, the  manifold $(2k+1)\mathbb{C}P^2\#l\overline{\mathbb{C}P^2}$ with $k>0$,  $l\geq0$ and  any almost complex structure satisfies the relations $$b_1=2h^{1,0}=2\tilde{h}^{1,0}=0\mbox{ and  } \tilde h^{0,1}>0,$$
hence violating all the equalities in \eqref{eqn-b1-hodge-complex-surface}.

\end{document}